\documentclass[]{elsarticle}

\usepackage{amssymb}
\usepackage{booktabs}
\usepackage{graphicx}
\usepackage[table,dvipsnames]{xcolor}
\usepackage{hyperref}
\usepackage{amsfonts}
\usepackage{amsmath}
\usepackage{mathtools}
\usepackage{graphicx}
\usepackage{xcolor, colortbl}
\usepackage{stmaryrd} 
\SetSymbolFont{stmry}{bold}{U}{stmry}{m}{n} 
\usepackage{algpseudocode}
\usepackage{algorithm2e}
\usepackage{enumitem}
\usepackage{multirow}
\usepackage{adjustbox}

\usepackage{subcaption}

\usepackage{tikz}
\usepackage{tikz-cd}

\newcommand{\R}{\mathbb{R}}
\newcommand{\Z}{\mathbb{Z}}

\newcommand{\kk}{\kappa}

\newcommand{\cM}{\mathcal{M}}

\newcommand{\cO}{\mathcal{O}}

\newcommand{\cMXZ}{ \cM_{\scst X}^{\scst Z}}
\newcommand{\Df}{D}

\newcommand{\dZ}{d^{\scst Z}}

\newcommand{\scst}{\scriptscriptstyle}
\newcommand{\Sf}{S^{\scst \Df}}
\newcommand{\mf}{m^{\scst \Df}}
\newcommand{\mU}{m^{\scst U}}
\newcommand{\mV}{m^{\scst V}}

\newcommand{\kkI}{\kappa_{\scst I}}

\DeclareMathOperator{\B}{B}
\DeclareMathOperator{\Ho}{H}
\DeclareMathOperator{\PH}{PH}

\DeclareMathOperator{\VR}{VR}

\DeclareMathOperator{\id}{id}

\DeclareMathOperator{\Rep}{Rep}
\DeclareMathOperator{\TMT}{TMT}
\DeclareMathOperator{\MST}{MST}

\newcommand{\rhoV}{\rho^{\scst V}}
\newcommand{\rhoU}{\rho^{\scst U}}

 \newtheorem{theorem}{Theorem}[section]
 
 \newtheorem{example}[theorem]{Example}
  \newtheorem{proposition}[theorem]{Proposition}
 \newdefinition{definition}[theorem]{Definition}
  \newdefinition{remark}[theorem]{Remark}

 \newproof{proof}{Proof}

\begin{document}

\begin{frontmatter}
\title{Topological Quality of Subsets via Persistence Matching Diagrams}
\author{Álvaro Torras-Casas\fnref{fn1,a1}}
\ead{atorras@us.es}
\author{Eduardo Paluzo-Hidalgo\fnref{fn1,a2,a1}}
\ead{epaluzo@\{uloyola, us\}.es}
\author{Rocio Gonzalez-Diaz\corref{cor1}\fnref{fn1,a1}}
\ead{rogodi@us.es}
 \cortext[cor1]{Corresponding author}
 \fntext[fn1]{A. Torras-Casas ORCID iD: 0000-0001-9937-0033.
 E. Paluzo-Hidalgo ORCID iD: 0000-0002-4280-5945.
 R. Gonzalez-Diaz ORCID iD: 0000-0002-5099-6294.}
\affiliation[a1]{organization={Universidad de Sevilla},
            addressline={Campus Reina Mercedes}, 
            postcode={41012},
            country={Spain}}
\affiliation[a2]{organization={Universidad Loyola Andalucia},
            addressline={Campus Sevilla}, 
            postcode={41704},
            country={Spain}}
\begin{abstract}
Data quality is crucial for the successful training, generalization and performance of machine learning models.
We propose to measure 
the quality of a subset concerning the dataset it represents,
using topological data analysis techniques.
Specifically, we 
define the persistence matching diagram, a topological invariant derived from combining embeddings with
persistent homology.
We provide an algorithm to compute it using minimum spanning trees.
Also, the invariant allows us to understand whether the subset ``represents well" the clusters from the larger dataset or not, and we also use it to estimate bounds for the Hausdorff distance between the subset and the complete dataset. 
In particular, 
this approach enables us to explain why the chosen subset
is likely to result in
poor performance
of a supervised learning model. 
\end{abstract}

\end{frontmatter}


\section{Introduction}

Data-driven machine learning (ML) is a subfield of artificial intelligence (AI) where models learn patterns towards making predictions and deriving insights from data. Data collection and preparation, model building, and model evaluations are the main parts of the workflow when dealing with data-driven models~\cite{9705125}.
Typically, about 80\% of the dataset is taken for training while reserving the remaining 20\% for testing (or a similar proportion).
However, the choice of the training data and validation are crucial~\cite{9857128} and usually forgotten.
Often, many training instances are very similar, and the training dataset could be taken in a smaller proportion; i.e. the model is ``learning'' what it already knows and part of the data could be removed as shown in~\cite{Li2023} in the case of large material datasets or studied in general in~\cite{10.12688/openreseurope.17554.1}. 
In addition, large training datasets could bring risks and environmental impacts in some areas, such as natural language processing. 
See, for example,~\cite{parrots}. 

Regarding the quality of training data,  researchers have observed that a small ``well-chosen"  dataset can be more effective than relying on a much larger but biased dataset \cite{upv}. In fact, including all available data in model training can occasionally lead to worse performance. See, for example,  \cite{10.1145/3592616} for a review of data quality requirements in ML development pipelines. 
Alternatively, towards trust in AI methods, the field of explainable AI focuses on different aspects of AI models such as: 
(1) inner explainability of the model (for example~\cite{alvarez2018senn}); 
(2) trustable modelling of the problem, for example~\cite{OLIVASPADILLA2024110418}; 
(3) post-hoc explainability methods such as LIME~\cite{10.1145/2939672.2939778} or SHAP~\cite{10.5555/3295222.3295230}; 
and (4) intrinsic explainable models such as classic decision trees and new variants such as, for example, \cite{LABER2023109239} or novel inventions on black-box models such as in \cite{WEI2024109991}. 
Besides, explainability methods vary significantly
depending on what they aim to explain and the features they used to provide those explanations. See~\cite{BAI2021108102} for explainable AI techniques used for pattern recognition
and~\cite {A2023100230} for a review on explainable AI models across various applications, not just pattern recognition.

Contrary to the usual explainable approaches, our study aims explainability directly from data and in a pre-trained phase, trying to identify weaknesses in the training dataset that could lead to bias with respect to the full dataset.
Specifically, in this paper, we present a tool that evaluates the \textit{topological} representativeness of a subset with respect to the full dataset and show that this helps analyze the quality of training datasets.
For this, we introduce the concept of persistence matching diagram,
which relies on block functions from~\cite{matchings} built on top of persistent homology. 
Such block functions are based on the idea of matching persistent homology barcodes, which has been applied to topological bootstrap~\cite{Reani2023} or to image segmentation~\cite{Stucki2023}.

Regarding the use of topological data analysis (TDA) to assess the quality of datasets in ML, 
\cite{Tudoreanu2022}
uses TDA to automatically detect data quality faults in datasets, addressing potential issues that might impact ML performance. 
The author proposed to 
transform a dataset into a multidimensional point cloud of real numbers that is based on calculating the distance between a subset of selected records in the dataset and all other records and then analyze the  Morse-Smale complexes of the point cloud depending on a parameter $\beta$.

Our study considers a pair of finite metric spaces $X$ and $Z$ together with an inclusion $X\subseteq Z$. 
In computational topology, the usual procedure to study the similarities between datasets starts by the computation of persistent homologies $\PH_*(X)$ and $\PH_*(Z)$, obtained from the respective Vietoris-Rips filtration, which are stable invariants~\cite{oudot}.
More precisely, one computes the interval decompositions (also called barcodes), $\B(\PH_*(X))$ and $\B(\PH_*(Z))$, and then, one compares them by some distance, such as the bottleneck or Wasserstein distance (see, for example~\cite{edelsbrunner} for an introduction to computational topology, focusing on the theoretical foundations and algorithms, and~\cite{carlsson} for TDA applications in data science and machine learning).
However, such comparisons might differ substantially from the underlying distribution of the datasets $X$ and $Z$,
since both distances are based only on combinatorial comparisons between the interval decompositions of $\PH_*(X)$ and $\PH_*(Z)$. 
This is the reason why, in this paper, we propose a different approach consisting of an indicator of topological data quality based on the block function $\cMXZ$ introduced in~\cite{matchings}, where there is an inclusion $X\subseteq Z$. 
Besides, we focus on 0-dimensional persistent homology of Vietoris-Rips filtrations, since this case has less technical difficulties than taking higher dimensions or other filtrations. 
In this paper, we show that we can compute $\cMXZ$ using minimum spanning trees, leading to an efficient and simple method.
Also, we use $\cMXZ$ to detect when a subset ``captures well" the connected components of a given dataset. 
In addition, we prove that we can use $\cMXZ$ to obtain the bounds of the Hausdorff distance between $X$ and $Z$. 
This further supports the consideration of $\cMXZ$ as a measure of the topological representativity of a dataset, offering greater flexibility compared to the Hausdorff distance. 

We start by reviewing the concept of block functions in Section~\ref{sec:matchings}, where we include examples and a novel description of how to compute it in the $0$-di\-men\-sion\-al case.
In Section~ \ref{sec:algorithm}, we provide a pseudocode for the computation of the 0-dimensional induced block function $\cMXZ$.
In Section~\ref{sec:tq}, we study the topological data quality of a subset with respect to the given dataset and give bounds for the Hausdorff distance between the subset and the full dataset. 
Readers interested in applying the proposed tool to ML problems can go directly to Section~\ref{applications}, where it is demonstrated through two different experiments. 
Finally, Section~\ref{sec:future} ends the paper with conclusions and future work.


\section{Block functions between barcodes induced by 
inclusion maps}
\label{sec:matchings}

Consider a finite metric space $Z$ and a subset $X$ of $Z$. 
Since distances are defined between the samples from $Z$, it is natural to consider the Vietoris-Rips filtration $\VR(X)$ and $\VR(Z)$ to create combinatorial models on the finite metric space. 
Let us remark that, in real applications, $Z$ would be a dataset and $X$ a subset of it, for example, the so-called training dataset. 
However, we stick to that mathematical notation from a topological context and keep a more general framework in the paper until Section~\ref{applications}.
Since we focus on the 0-dimensional persistent homology, we consider the 1-skeleton of Vietoris-Rips filtration (also denoted as $\VR$ to unload the notation) and fix $\Z_2$ as the ground field.
Throughout this text, we use the notation $U=\PH_0(X)$ and $V=\PH_0(Z)$ for the 0-dimensional persistent homology over $\Z_2$ of $\VR(X)$ and $\VR(Z)$ respectively. 
In this section,
our goal is to explain how to compute the block function $\cMXZ$ from the barcode $\B(V)$ to the barcode $\B(U)$ induced by the inclusion $X\subseteq Z$.


\subsection{Background}

First, we introduce some technical background from the computational topology field.
All the information can be found, for example,  in \cite{edelsbrunner, OudotBook2015,Chazal2021}.


\subsubsection{Finite metric spaces and the Hausdorff distance}

Let $M$ be a finite metric space with metric $d^{\scst M}$.  
Given a pair of subsets $A,B\subseteq M$, we define their pairwise distance as the quantity 
\[
d^{\scst M}(A,B) = \min\big\{d^{\scst M}(a,b) \mbox{ for all } a \in A \mbox{ and } b \in B\big\}\,.
\] 
Given $A \subseteq M$ and a point $x \in M$, we use the notation $d^{\scst M}(x, A) = d^{\scst M}(\{x\}, A)$. 

A disadvantage of the definition of distance $d^{\scst M}$ between subsets is that it is not effective in distinguishing subsets that are potentially very different since, for example, if $A\cap B\neq \emptyset$ then $d^{\scst M}(A,B)=0$. 
For this purpose, it is more useful to consider
the \emph{Hausdorff} distance
\[
d^{\scst M}_{\scst H}(A,B) = \max\big\{\max\big\{ 
d^{\scst M}(a, B) \, \vert \,
a \in A
\big\}, \max\big\{ 
d^{\scst M}(A, b) \, \vert \,
b \in B
\big\}\big\}\,,
\]
since $d^{\scst M}_{\scst H}(A,B)=0$ if and only if $A=B$.


\subsection{The 1-skeleton of the Vietoris-Rips filtration, $\VR(Z)$}

Given a dataset $Z$ with metric $\dZ$,  the 1-skeleton $\VR(Z)$ of the Vietoris-Rips filtration is a family of graphs 
$\VR(Z)=\big\{\VR_r(Z)\big\}_{r\in[0,\infty)}\,,
$
being $r$ a scale parameter, satisfying that there is an injective morphism of graphs $\VR_r(Z)\xhookrightarrow{} \VR_{s}(Z)$ for $r\leq s$.
Specifically, fixed $r\geq 0$, $\VR_r(Z)$ is a graph with vertex set $Z$ containing the edge $[x,y]$ with endpoints $x,y\in Z$, whenever 
the distance between $x$ and $y$ is less than or equal to $r$, that is, $\dZ(x,y)\leq r$.

\begin{example}\label{ex:VR-complex-X-Z}
    Consider Figure~\ref{ex:0-homology-VR} (top row) where a 
   set $Z$ consisting of six points, with a sample $X$ with three red points, is pictured. 
    As we increase the filtration value $r \in \{0.0, 1.0, 1.1, 1.2, 2.5\}$, $\VR_r(X)$ and $\VR_r(Z)$ increase. 
    We indicate the edges of $\VR_r(X)$ in red and the remaining edges in blue. 
\end{example}


\subsubsection{The 0-dimensional homology group of  $\VR_r(Z)$}\label{subsec:0-homology}

Fixed $r\geq 0$, the 0-dimensional homology group of $\VR_r(Z)$, is a quotient group, denoted as  $\Ho_0(\VR_r(Z))$ that consists of the vector space generated by the classes associated with the connected components of $\VR_r(Z)$.
Specifically, 
$\Ho_0(\VR_r(Z))=\big\langle \pi_0(\VR_r(Z))
\big\rangle$ is the free $\Z_2$-vector space generated by the set $\pi_0(\VR_r(Z))$ of connected components of $\VR_r(Z)$. 

To describe $\pi_0(\VR_r(Z))$, suppose that $Z=\{z_0,\dots,$ $z_{n-1}\}$.  
Each $\alpha_r \in \pi_0(\VR_r(Z))$ is associated with a set of points $A_{r}\subseteq Z$ such that:

\begin{itemize}
\item [1)] for any two points $x,y\in A_{r}$, there exists a \emph{path} $\{x_0,x_1,\dots,x_{\ell}\}\subseteq A_{r}$ satisfying that
$x_0=x$, $x_{\ell}=y$ and 
$\dZ(x_{i},x_{i+1})\leq r$ for all 
$i\in \llbracket \ell\rrbracket$,
where the expression
 $\llbracket \ell\rrbracket$ denotes the set $\{0,1,\dots,\ell-1\}$; and 
\item[2)] $d^{\scst Z}(z,A_r)>r$ for any point $z\in Z\setminus A_r$.
\end{itemize}

Besides, given $z\in Z$, we denote by $A_{r}(z)$ the set of points of $Z$ that are connected to $z$ in $\VR_r(Z)$. 
From this viewpoint, $\pi_0(\VR_r(Z)) = Z/\! \sim$ where $\sim$ is the equivalence relation given by $x  \!\sim \!   y$ if and only if $x \in A_r(y)$.
Finally,  
$\alpha_r\in \pi_0(\VR_r(Z))$ associated with the connected component $A_r$ is always represented by the point of $A_r$ with the lowest label. 
That is, $\alpha_r=[z_j]$ for some 
$z_j\in A_r$ such that 
$j=\min\big\{\ell\in \llbracket n\rrbracket
\ |\ z_{\ell}\in A_r\big\}$.

\begin{figure}[ht!]
    \centering
    \includegraphics[width=0.8\textwidth]{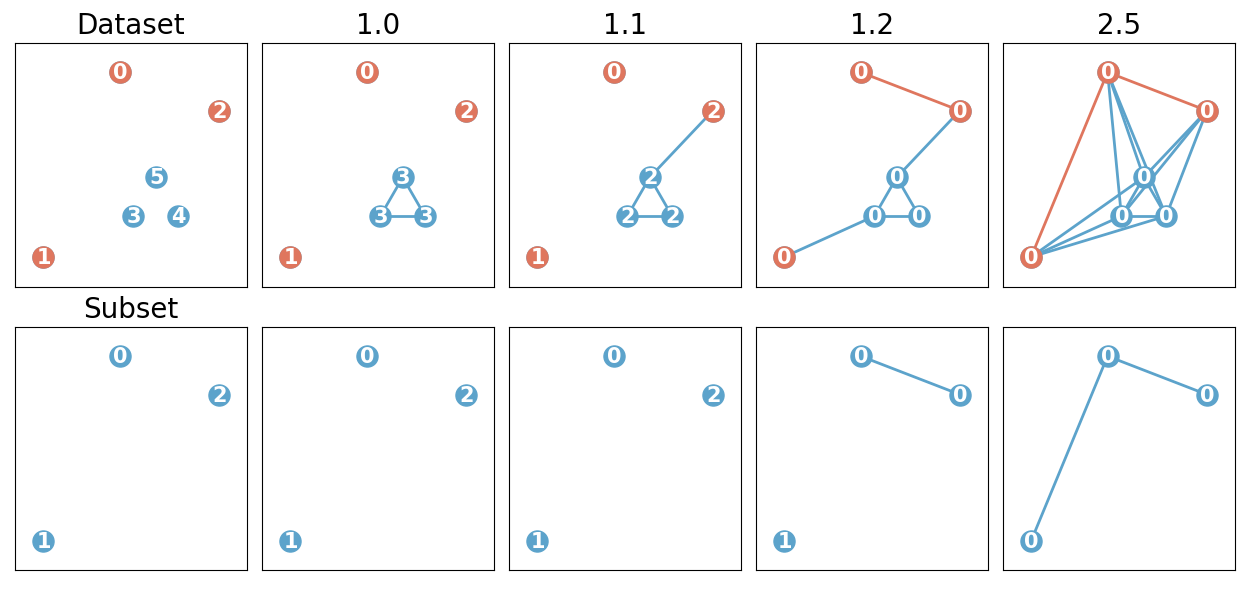}
    \caption{Vietoris-Rips filtration of a set $Z$ (top row) and a subset $X$ (bottom row) with connected components indexed by non-negative integers.}
    \label{fig:0-homology-VR}
\end{figure}

\begin{example}\label{ex:0-homology-VR}
    Consider the same pair $X\subset Z$ from Example~\ref{ex:VR-complex-X-Z}, and suppose that we have labeled the points from $X$ and $Z$ as depicted on the left of Figure~\ref{fig:0-homology-VR}. 
    Thus, we might write $Z=\{z_0,z_1,z_2,z_3,z_4,z_5\}$ and $X=\{z_0,z_1,z_2\}$. 
    In the same Figure~\ref{fig:0-homology-VR}, we consider filtration values ranging over $r \in \{1.0,1.1,1.2, 2.5\}$ and we label the vertices of $\VR_r(Z)$ (top row) and $\VR_r(X)$ (bottom row) with the smaller index in their respective connected components. 
    For example, we have 
    $
    \Ho_0(\VR_{1.2}(Z))=\big\langle[z_0]\big\rangle
    $ 
    while 
    $
    \Ho_0(\VR_{1.2}(X))=\big\langle[z_0],[z_1]\big\rangle
    $.
\end{example}


\subsubsection{The 0-dimensional persistent homology of $\VR(Z)$ and merge trees}\label{subsec:pershom}

Persistence modules are an algebraic generalization for persistent homology.
Specifically, a persistence module $V$ indexed by $\R$ consists of a set of vector spaces $\big\{V_t\big\}_{t\in\R}$ and a set of linear maps 
$\big\{\rhoV_{st}\colon V_s \rightarrow V_t\big\}_{s\leq t}$, called the \emph{structure maps} of $V$,  satisfying that $\rhoV_{jt}\rhoV_{ij} = \rhoV_{it}$ and $\rhoV_{tt}$ being the identity map, for $i\leq j\leq t\in\R$.

Given an interval $I=[a,b)\subset \R$,
the \emph{interval module}, $\kkI$, is a particular case of persistence module consisting of $(\kkI)_t= \Z_2$ for all $t \in I$ and $(\kkI)_t= 0$ otherwise, while the linear map $\rho^{\kkI}_{ij}\colon (\kkI)_i\to (\kkI)_j$ is the identity map whenever $a\leq i\leq j< b$. 
In this work, for most cases, we consider only interval modules $\kkI$ for $I=[a,b)$ such that $a=0$.
Thus, for ease, we often denote an interval module $\kk_{[0,b)}$ by the simpler notation $\kk_{b}$ for all $b>0$.
Also, we denote by  $\kk_{\infty}$  
the persistence module that consists of 
$(\kk_{\infty})_t= \Z_2$ for all $t\geq 0$ and 
$(\kk_{\infty})_t=0$ otherwise; and $\rho^{\kk_\infty}_{ij}\colon(\kk_{\infty})_i\to (\kk_{\infty})_j$ being the identity map whenever $0\leq i\leq j$ and the zero map otherwise.

The 0-dimensional persistent homology of $\VR(Z)$, denoted as $U=\PH_0(Z)$, 
\emph{encapsulates} the 0-dimensional topological events that occur within the filtration. 
Specifically, $U$ is a \emph{persistence module} given by the set  of 0-dimensional homology groups 
$\big\{U_r=\Ho_0(\VR_r(Z))\big\}_{r\in[0,\infty)}$
and the set of linear maps 
\[
\big\{\rhoU_{rs}\colon\Ho_0(\VR_r(Z))\to \Ho_0(\VR_s(Z))\big\}_{r\leq s}
\]
that are induced by the inclusion map
$\VR_r(Z)\xhookrightarrow{} \VR_s(Z)$. 

In particular, given $\alpha_r=[z_j]
\in \pi_0(\VR_r(Z))$, we have that:
 \[\rhoU_{rs}([z_j])=
 \left\{
 \begin{array}{cl}
 \mbox{$[z_j]$} &
 \mbox{  if $\dZ(A_s(z_j),A_s(z_{\ell}))> s$ for all  $\ell\in [\![j]\!]$,}
 \\
  \mbox{$[z_{i}]$}  &
 \mbox{  for some  $i\in [\![j]\!]$ such that $A_s(z_j)=A_s(z_i)$.}
 \end{array}
 \right.
 \]
Thus, if $\rhoU_{rs}([z_j])=[z_j]$ then 
$A_s(z_j)$ is represented by $z_j$. 
On the other hand, if $\rhoU_{rs}([z_j])=[z_i]$ for some $i\in [\![j]\!]$ then  $A_s(z_j)$ is represented by $z_i$.

The connected components in $\VR(Z)$ start being the isolated points of the whole dataset $Z$ and, as the filtration parameter increases, $\PH_0(Z)$ records the death values of such components. 
In this way,  all classes in $\PH_0(Z)$ are \emph{born} at $0$, meaning that 
$
\pi_0(\VR_0(Z))=
\big\{\ [z_i]\ |\ i\in \llbracket n\rrbracket\ \big\}
$.
On the other hand, a class $[z_j]\in \Ho_0(\VR_0(Z))$ is said to \emph{die} at $b>0$ 
if: 
\begin{itemize}
    \item[1)] $d^{\scst Z}(A_r(z_j),A_r(z_{\ell}))>r$ for all  $\ell\in [\![j]\!]$ and $r \in [0,b)$; and 
\item[2)] $A_b(z_j)=A_b(z_i)$ for some $i\in [\![j]\!]$. 
\end{itemize}
Observe that if $A_b(z_j)=A_b(z_{i})$ for some $i\in [\![j]\!]$ then
$
\rhoU_{0b}([z_j])=\rhoU_{0b}([z_{i}])=[z_i]
$, concluding  that the class $\alpha_0=[z_j]+[z_{i}]\in \Ho_0(\VR_0(Z))$ belongs to $\ker\rhoU_{0b}$. 

We then 
 consider the set of triplets $(z_j, b_j, z_i) \in Z\times \R_{>0} \times Z$ such that
 $[z_j]\in \Ho_0(\VR_0(Z))$ dies at value $b_j>0$ and $\rhoU_{0b_j}([z_j]) = [z_i]$. Here we ignore the component $[z_0]$ which never \emph{dies}. 
We denote such set as $\TMT(Z)$ and call it the \emph{triplet merge tree} for $Z$.
Notice that this is a variant of the original definition from~\cite{triplets}, which we modify for convenience in our case as there is no total order in the edge set, and several edges might appear at the same filtration value.

\begin{figure}[ht!]
    \centering    \includegraphics[width=0.7\textwidth]{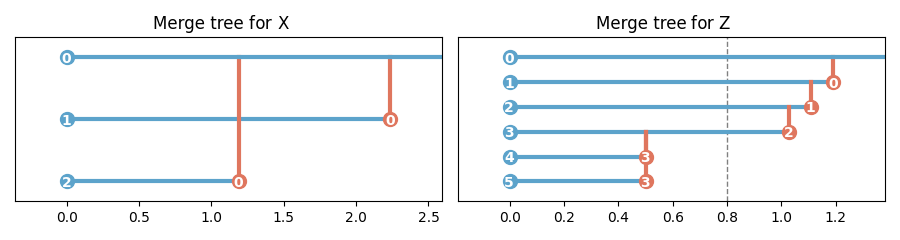}
    \caption{Merge tree representation for $X$ and $Z$. 
    Triplets $(z_j,b_j,z_i)$ are plotted as a blue horizontal segment labeled by $j$ on the left and $i$ on the right, followed by a red vertical segment 
    that connects to the blue horizontal segments labeled on the right by $i$. 
    In addition, on top of the representation, we plot a horizontal blue line corresponding to the component $[z_0]$, which never dies.     
    The blue horizontal intervals compound the barcodes of $\PH_0(X)$ (merge tree representation on the left) and $\PH_0(Z)$ (merge tree representation on the right).     
    }
    \label{fig:merge-trees-X-Z}
\end{figure}

\begin{example}\label{ex:pershom-merge-trees}
    We continue Example~\ref{ex:0-homology-VR} and track the connected components from both $\VR(X)$ and $\VR(Z)$. 
    For example,  $\TMT(X)$ contains the triplets $(z_2, 1.19, z_0)$ and $(z_1, 2.24, z_0)$.
    Such triplets are depicted on the left of  Figure~\ref{fig:merge-trees-X-Z} such that, for example,  triplet $(z_2, 1.19, z_0)$ is plotted as a blue horizontal segment labeled by $2$ on the left and $0$ on the right, followed by a red vertical segment that connects to the horizontal segment of $0$.  
    The resulting representation encodes the merge tree information from $\VR(X)$.
    On the right of Figure~\ref{fig:merge-trees-X-Z}, we repeat the same principle using $Z$. 
    An advantage of representing merge trees is that we know instantly $\Ho_0(\VR_r(Z))$ for all $r\geq 0$. 
    For example, choosing $r = 0.8$, we consider the dashed vertical line on the 
    merge tree representation pictured on the right of Figure~\ref{fig:merge-trees-X-Z}; such line crosses four blue segments so that $\Ho_0(\VR_r(Z))=\big\langle
    [z_0],[z_1],[z_2],[z_3]\big\rangle
    $.
\end{example}

One reason to consider triplet merge trees is that, using a minimum spanning tree for $\VR(Z)$, one can obtain $\TMT(Z)$ directly, as we explain briefly in Section~\ref{sec:algorithm}.
Another reason is that they help understand the operators $\ker^\pm_{b}$ that will be used in Subsection~\ref{sec:block-function-computation}. 
Basically, given the persistence module $U=\PH_0(Z)$ with structure maps $\rho^{\scst U}$, 
the operators $\ker^\pm_{b}$ are defined, 
for all $b> 0$, as
$\ker^-_b(U
)= \bigcup_{0 \leq r < b} \ker(\rho^{\scst U}_{0r})$
and 
$\ker^+_b(U)
= \ker(\rho^{\scst U}_{0b})
$,
and encapsulate the classes in $\PH_0(Z)$ that die at $r$ for $r<b$ and  $r\leq b$, respectively.
Besides, observe that
$ \ker^-_b(U
)\subseteq \ker^+_b(U)\subseteq U_0\,.$
Using $\TMT(Z)$, we might write these operators as follows
\[
\ker^-_b(U) = \big\langle \; [z_i]+[z_j] \;\big\vert \;(z_j,b_j,z_i) \in \TMT(Z) \mbox{ and } 
b_j< b \;\big\rangle
\]
and 
\[
\ker^+_b(U) = \big\langle\; [z_i]+[z_j]\; \big\vert \; 
(z_j,b_j,z_i) \in \TMT(Z) \mbox{ and } 
b_j\leq b \; \big\rangle\,.
\]

\begin{example}
    Recall, from Example~\ref{ex:pershom-merge-trees}, that $\TMT(X)$ consists of triplets $(z_2, 1.19, z_0)$ and $(z_1, 2.24, z_0)$, depicted on the left of Figure~\ref{fig:merge-trees-X-Z}.
    Then, for $V=\PH_0(X)$,  we have:
    \[
    \ker^-_{2.24}(V) = \langle [z_2] + [z_0]\rangle
    \;\;\mbox{ and }\;\;
    \ker^+_{2.24}(V) = \langle [z_2] + [z_0], [z_1]+[z_0]\rangle\,.
    \]
\end{example}


\subsubsection{Barcodes, persistence diagrams and multisets
}

It is known that, in general, a persistence module $V$ has a unique decomposition as a direct sum of interval modules (see~\cite{crawley-boevey-2015}). 
The interval modules are in bijection with intervals over $\R$, so 
$V$ is uniquely characterized by a \emph{multiset} called barcode.

A multiset is a pair $(S,m)$ composed of a set $S$ together with an assignment $m\colon S\rightarrow \Z_{>0} \cup \{\infty\}
$ that maps elements from $S$ to their multiplicity. 
The representation of a multiset $(S,m)$ is the set 
\[\Rep\, (S,m) = \big\{ \; (s,\ell) \;\vert\; s \in S
\mbox{ and } 
\ell\in \llbracket m(s)\rrbracket\; \big\}\,.\]
In particular,
the cardinality of $(S,m)$ is
$\# (S,m) = \# \Rep\, (S,m)$.

The persistence modules we 
deal with in this paper are uniquely characterized by multisets of intervals that involve finite sets and finite multiplicities. 
Thus, given $U=\PH_0(Z)$, there is a multiset $\B(U)=(S^{\scst U},\mU)$,
where $S^{\scst U}$ is a set of intervals over
$\R$ together with an assignment called \emph{multiplicity},
$\mU:S^{\scst U}\rightarrow \Z_{> 0}$, satisfying that there is
a persistence isomorphism 

\[
\mbox{$
U\;
\simeq\; \bigg(\bigoplus_{\scst  (J,\ell) \in \Rep\B(U)}  \kk_{\scst J}\bigg) \oplus \kk_\infty
\;=\;
\bigg(\bigoplus_{\scst J \in S^{\scst U}} 
\bigoplus_{\scst \ell\in \llbracket \mU\!(J)\rrbracket}
\kk_{\scst J}\bigg)\oplus \kk_\infty$.}
\]
Specifically, its barcode, $\B(U)=(S^{\scst U},\mU)$, is such that all intervals from $S^{\scst U}$ are of the form $[0,b)$ for values $b>0$. 
Notice that, in this paper, we do not consider the infinity interval $[0,\infty)$ as an element from $\B(U)$.
This is why we might consider $S^{\scst U}$ as a subset of $\R$, where an interval $[0,b)$ is substituted by its endpoint $b\in \R_{>0}$. 
Similarly, we also write $\mU\!(b)$ and $\kk_b$ instead of $\mU\!(J)$ and $\kk_{\scst J}$ without ambiguities.
Besides, since $Z$ is finite, we have that $\mU\!(b)<\infty$.
Then, 
\begin{equation}\label{eq:decomposition-U-mUb}
\mbox{$
U\;\simeq\; 
\bigg(\bigoplus_{\scst b>0} 
\bigoplus_{\scst \ell \in \llbracket \mU\!(b) \rrbracket}
\kk_b \bigg) \oplus \kk_\infty
$}.
\end{equation}
These intervals track when connected components merge.
In particular, there is a bijection between $\Rep\B(U)$ and $Z\setminus \{z_0\}$ in the sense that, for each $b>0$ and  
$\ell \in \llbracket \mU\!(b) \rrbracket$, the pair $([0,b),\ell)$ can be univocally associated with  a class $[z_j]\in\PH_0(Z)$, for $j\in[\![n]\!]$, that dies at $b$. 

\begin{remark}\label{rem:iso-TMT-RepU}
    There is an isomorphism $\TMT(Z)\simeq \Rep \B(U)$ given by sending a triplet $(z_j, b_j, z_i)\in \TMT(Z)$ to 
    $([0, b_j), k) \in \Rep \B(U)$ for some $k \in \llbracket \mU\!(b_j) \rrbracket$. 
\end{remark}

\begin{example}\label{ex:barcode-PH0-Z}
    Here we consider $U=\PH_0(Z)$ from Example~\ref{ex:pershom-merge-trees} again. 
    Then,
    \[\TMT(Z)=\big\{(z_5,0.50,z_3), (z_4,0.50,z_3), (z_3, 1.03, z_2), (z_2, 1.11, z_1), (z_1, 1.19, z_0)\big\}.\] 
    In other words, both $[z_5], [z_4] \in U_0$ die at value $0.50$,  $[z_3] \in U_0$ dies at value $1.03$, $[z_2] \in U_0$ dies at value $1.11$ and $[z_1] \in U_0$ dies at value $1.19$; 
    while the class $[z_0] \in U_0$ never dies. 
    Thus, 
    \[
    U \simeq 
    \kk_{0.50}\oplus \kk_{0.50}    
    \oplus \kk_{1.03} \oplus \kk_{1.11} \oplus \kk_{1.19} \oplus \kk_\infty\,.
    \] 
    That is, $\B(U)=(S^{\scst U}, \mU)$ with $S^{\scst U}=\{0.50, 1.03, 1.11, 1.19\}$ and also $\mU\!(0.50)$ $=2$ and $\mU\!(1.03)=\mU\!(1.11)=\mU\!(1.19)=1$. 
    A representation of the barcode $\B(U)$ is plotted by the finite horizontal blue segments on the right in Figure~\ref{fig:merge-trees-X-Z}.
    \end{example}
    Barcodes can also be represented using \emph{persistence diagrams}, consisting of a set of points $(p,q)$ in the cartesian plane representing a homology class born at $p$ which dies at $q$. 
        
\subsection{From inclusion maps to block functions: the induced block function $\cMXZ$
}~\label{sec:block-function-computation}

In~\cite{matchings}, the authors defined a block function
$
\cMXZ:S^{\scst V} \times S^{\scst U} \rightarrow \Z_{\geq 0}
$ 
induced by a given persistence morphism $f:V\to U$ and provided a matrix-reduction algorithm to compute it. 
In our present case, since we work with 0-dimensional persistence homology and persistence morphisms induced by inclusions, the definition of the block function is simplified as follows. 
For a pair of sets $X\subseteq Z$, and  pair of intervals $I=[0,a)$ and $J=[0,b)$ we have
\begin{equation}\label{equation:bf-formula}
\cMXZ(I,J) = 
\dim\bigg(
\dfrac{
\ker_{a}^+(V) \cap \ker_{b}^+(U)
}{
\ker_{a}^-(V) \cap \ker_{b}^+(U) + \ker_{a}^+(V) \cap \ker_{b}^-(U)
}
\bigg)\,,
\end{equation}
where we used the fact that $\ker^-_a(V), \ker^+_a(V)\subseteq V_0$ together with $V_0\subseteq U_0$.
Given $a,b>0$, we write $\cMXZ(a,b)$ instead of $\cMXZ([0,a), [0,b))$ to simplify notation.

Observe 
that, given $a<b$, we have $\ker^+_a(V) \subseteq \ker^-_b(U)$ and this implies that $\cMXZ(a,b)=0$.
Besides, $\cMXZ$ is well-defined as 
a consequence of the following remark obtained combining Remark~\ref{rem:iso-TMT-RepU} 
that states that $\TMT(X)\simeq \Rep \B(V)$ and $\TMT(Z)\simeq \Rep \B(U)$
with Proposition~\ref{prop:Mf-bijection-TMTX-TMTXZ} knowing that $\TMT(X,Z)\subseteq \TMT(Z)$. Composing both isomorphisms and the inclusion we obtain an injection $\Rep \B(V)\hookrightarrow \Rep \B(U)$ sending $\cM^Z_X(a,b)$ copies of $[0,a)$ to copies of $[0,b)$. Consequently, we have,
\begin{remark}~\label{rem:matching-well-defined}
    $\sum_{b' \leq a} \cMXZ(a,b') = \mV(a)$ and $\sum_{b \leq a'} \cMXZ(a',b) \leq \mU(b)$.
\end{remark}

An algorithm for computing $\cMXZ$ induced by an inclusion $X \subseteq Z$ is given in Section~\ref{sec:algorithm}. 
Let us see now an example.

\begin{example}\label{ex:matching-example}
    We consider $X\subseteq Z$ from Example~\ref{ex:0-homology-VR} and recall the computation of $\B(U)$ from Example~\ref{ex:barcode-PH0-Z}, that is, $\B(U)=(S^{\scst U}, \mU)$ with $S^{\scst U}=\{0.50, 1.03, 1.11, 1.19\}$ and  $\mU\!(1.03)=\mU\!(1.11)=\mU\!(1.19)=1$  and $\mU\!(0.50)$ $=2$.
    Also, from Example~\ref{ex:pershom-merge-trees}, we can observe that $\B(V)=(S^{\scst V}, \mV)$ is such that $S^{\scst V}=\{1.19, 2.24\}$ with multiplicities $\mV(1.19)=\mV(2.24)=1$. 
     Using Formula~(\ref{equation:bf-formula}),
    $\cMXZ$ is nonzero only for 
    $
    \cMXZ(1.19, 1.19)=\cMXZ(2.24,1.11)=1\,.
    $
\end{example}


\section{Matching diagrams for topological quality of subsets}\label{sec:tq}

In this section, we introduce $\Df=D(X,Z)$ a diagram representation for $\cMXZ$ induced by $X\subseteq Z$, which is analogous to persistence diagrams but with some fundamental differences. 
Besides, we interpret the matching diagram $\Df$ giving an intuition of the relation between $X$ and $Z$, to know whether $X$ ``represents well'' the clusters from $\VR_r(Z)$ for $r\geq 0$. 
Additionally, using $\cMXZ$, we provide bounds 
to the Hausdorff distance between $X$ and $Z$.

We define the \emph{0-dimensional persistence matching diagram} $\Df$, also called \emph{matching diagram} for short,  to be the multiset $(\Sf, \mf)$ given by a set $\Sf \subset \overline{\R} \times \R$, where $\overline{\R}=\R\cup \{\infty\}$, together with the multiplicity function $\mf$ given by 
\[\begin{array}{l}
\mbox{
$\mf((a,b))=\cMXZ(a,b)$ for all pairs $(a,b) \in \Sf \cap \R^2 \;
$ and
}
\\
\mbox{
$\mf((\infty,b)) = \mU\!(b)-
\sum_{a < \infty} 
\cMXZ
((a,b))$
 for all $b>0$\,.
 }
 \end{array}\] 

Notice that $(a,b) \in \Sf$ 
implies $0<b\leq a$.
Besides, $\Df=\Df_{\scst O}\sqcup \Df_{\scst \infty}$ where $\Df_{\scst O}$ is the multiset  given by the pairs $(a,b) \in S^D$ with $a<\infty$ and with multiplicity $\mf((a,b))$; and $\Df_\infty$ is the multiset of pairs $(\infty,b) \in S^D$ with multiplicity $\mf((\infty,b))$. 

\begin{figure}[ht!]
    \centering
    \includegraphics[width=0.35\textwidth]{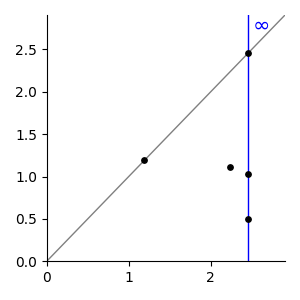}
    \caption{Depiction of the set $\Sf$ associated to the matching diagram $\Df$ detailed in Example~\ref{ex:diagram-matching}.}
    \label{fig:matching-diagram}
\end{figure}

\begin{example}\label{ex:diagram-matching}
    We consider $\cMXZ$ 
    from Example~\ref{ex:matching-example}.
    In this case, one can compute the matching diagram $\Df=(\Sf, \mf)$ and plot the points from $\Sf$ as done in Figure~\ref{fig:matching-diagram}. 
    Notice that the pairs $(\infty, b) \in \Sf$ have been plotted over a vertical blue line with the sign $\infty$ next to it. 
    Also, 
    $\mf((\infty, 0.50))=2$; all other points have multiplicity $1$.
\end{example}

\subsection{Interpreting $D(X,Z)$ to study the relation between $X$ and $Z$}\label{subsec:interpretation-Df}

In this subsection, we use $D=D(X,Z)$ to study the relationship between $X$ and $Z$.
First, given finite metric spaces $X\subseteq Z$, we assume that $Z$ has $n$ elements so that $X=\{z_0,\dots,z_{\ell-1}\}$ and $Z=\{z_0,\dots,z_{n-1}\}$ for some $\ell\leq n$.
We define $\TMT(X,Z)$ to be the subset of $\TMT(Z)$ given by triplets $(z_j,b,z_i)$ such that $i < j < \ell$; so that $z_i, z_j \in X$ by hypothesis.

To start, we relate $\#X$ and $\#Z$ with $\# \Df_{\scst O}$ and $\# \Df_{\infty}$. 
First, by Remark~\ref{rem:iso-TMT-RepU} and Remark~\ref{rem:matching-well-defined},
\begin{align*}    
\mbox{$
\#\TMT(X) 
= \sum_{[0,a)\in S^V} m^V(a) 
= \sum_{a \in \R} \Big( \sum_{b' \leq a} \cMXZ(b',a) \Big)$}\\
&\hspace{-3cm}\mbox{$= \sum_{(a,b) \in \R^2} \cMXZ(a,b)\,.$}
\end{align*}
Since $\TMT(X)$ contains a triplet for each point from $X$, except $z_0$, we have $\# X -1 = \#\TMT(X)$. 
Similarly, we obtain $\TMT(X, Z)= \# X -1$ and $\TMT(Z)=\# Z - 1$.
On the other hand, fixed $b>0$, the number of unmatched triplets $(z_i, b, z_j) \in \TMT(Z)\setminus \TMT(X,Z)$, by Remark~\ref{rem:matching-well-defined}, must be $\mU\!(b) -
\sum_{a < \infty} \cMXZ(a,b)$. Consequently, we have the equality
$$\# Z -\# X = \#\TMT(Z)  - \#\TMT(X,Z) 
= \mbox{$\sum_{b>0} \left( \mU\!(b) -
\sum_{a < \infty} \cMXZ(a,b)\right)\,.$}
$$
To summarize, using the definition of $\Df$, we obtain the following remark. 

\begin{remark}~\label{rem:dim-X-Z}
Given $X\subseteq Z$, consider the persistence morphism $f:V\to U$ for $V=\PH_0(X)$ and $U=\PH_0(Z)$. The following equalities hold:
\begin{itemize}
    \item $\# X -1= \sum_{(a,b)\in S^{\Df}\cap \R^2} \mf((a,b))\, =\# \Df_{\scst O}\,;$
    \item $\# Z-\# X = \sum_{(\infty, b) \in S^{\Df}} \mf((\infty,b)) =\# \Df_{\infty}\,.$
\end{itemize}
\end{remark}

Now, we focus on $\Df_{\infty}$ and consider the maximum $y$-coordinate of one of its points, that is, we take the 
quantity
\[
\eta_f = \max_{b\geq 0}\, \big\{ \, b \,\mid\, 
\mf((\infty,b)) >0 \mbox{ or } b=0\big\}\,.
\]
The following result quantifies how many connected components from  $\VR_r(Z)$  
are missed by $X$ and implies that the closer the points from $\Df_{\infty}$ are to $(\infty,0)$, the better.

\begin{proposition}\label{prop:components-infty-line}
Fixed $r\geq 0$, there are $
\sum_{r < b} m^{\Df}((\infty, b))
$
components in $\VR_r(Z)$ that contain no points from $X$. 
In particular, 
all components from $\VR_r(Z)$ contain points from $X$ for all $r \geq \eta_f$.
Besides, if $\eta_f=0$ then $X=Z$.
\end{proposition}

\begin{proof}
First, recall our indexing convention on $X$ and $Z$, so that triplets $(z_j, b, z_i) \in \TMT(Z)\setminus \TMT(X,Z)$ are such that $j\geq \ell$, where $\ell=\# X$. 
\\
In particular, given $r< b$, the triplet $(z_j, b, z_i) \in \TMT(Z)\setminus \TMT(X,Z)$ indicates that a component in $\VR_r(Z)$ is indexed by $z_j$, which is its vertex of minimum index
and $j\geq \ell$; 
so it is a component with no points from $X$. 
This implies that, for a fixed pair 
$r<b$, 
there are $\mf((\infty, b))$ connected components in $\VR_r(Z)$ with no points from $X$ and which correspond to triplets of the form $(\cdot, b, \cdot)$ in $\TMT(Z)\setminus \TMT(X,Z)$. 
Altogether, there are 
$
\sum_{r<b} \mf((\infty, b))
$
components in $\VR_r(Z)$ with no samples from $X$. 
\\
    If $\eta_f=0$ then, by Remark~\ref{rem:dim-X-Z}, $\# Z- \# X =0$ and since $X\subseteq Z$ then
     $X=Z$.
     \qed
\end{proof}

The following result indicates that the closer the points from $\Df_{\scst O}$ are to the diagonal axis, the more relevant the clusters from $X$ become within $Z$.

\begin{proposition}
    Given $s \geq 0$, suppose that $\cMXZ(a,b)=0$ for all $|b-a|> s$. 
    Then, for any $r\geq 0$, given an element $[z_i] \in \pi_0(\VR_r(X))$ such that $\rho_{r(r+s)}([z_i])=[z_i]$ then $f_r([z_i])=[z_i]$. 
    In particular, if $s=0$ then $f_r([z_i])=[z_i]$ for all $[z_i] \in \pi_0(\VR_r(X))$ and all $r\geq 0$.
\end{proposition}

\begin{proof}
    We prove this by showing the reciprocal statement, that is:
    given $[z_i] \in \pi_0(\VR_r(X))$ such that $f_r([z_i])=[z_j]$ for some $j<i$, then $\rho_{r(r+s)}([z_i])=[z_k]$ for some $k<i$. 
    By hypotheses, it follows that $[z_i]-[z_j] \in \ker(f_r)$, and we assume $[z_i]\neq [z_j]$ as elements in $\pi_0(\VR(X)_r)$, since otherwise the claim is direct. 
    Next, consider $[z_i]_0, [z_j]_0 \in \pi_0(\VR(X)_0)$ such that $\rho_{0r}([z_i]_0)=[z_i]$ and $\rho_{0r}([z_j]_0)=[z_j]$.
    There exist values $a,b>0$ such that $[z_i]_0-[z_j]_0 \in \ker^+_a(V) \setminus \ker^-_a(V)$ and $[z_i]_0-[z_j]_0 \in \ker^+_b(U) \setminus \ker^-_b(U)$; these values must exist by the definitions of $\ker^\pm_a(V)$ and $\ker^\pm_b(V)$.
    This implies that $\cMXZ(a,b) \neq 0$ and, by hypotheses, we have $a-b \leq s$. 
    Also, as $[z_i]-[z_j]\in \ker(f_r)$ and $[z_i]\neq [z_j]$, we have $b < r < a$
    and so $\rho_{r(r+s)}([z_i]-[z_j])=0$.
    Hence, there exists $k\leq j < i$ such that $\rho_{r(r+s)}([z_i])=\rho_{r(r+s)}([z_j])=[z_k]$. 
    \qed
\end{proof}

This way, we say that $X$ ``represents well" $Z$ if the points in $D_{\infty}$ are ``close" to $(\infty,0)$ and the points of $D_{\scst O}$ are ``close" to the diagonal axis. Furthermore:
\begin{itemize}
\item If $X=Z$ then all the points of $D_{\scst O}$ are in the diagonal axis and $D_{\scst \infty}=\emptyset$.
\item If $X$ is a \textit{$\varepsilon$-representative} subset of $Z$ (meaning that for any point $z\in Z$ there exists a point $x\in X$ such that $d^{\scst Z}(x,z)\leq \varepsilon$) then  any point $(a,b)\in D_{\scst O}$ satisfies that $|a-b|\leq 2\varepsilon$ and any point $(\infty,b)\in D_{\infty}$ satisfies that $|b|<2\varepsilon$.
 \end{itemize}
The second point is a consequence of applying the stability result of  block functions induced by embeddings~\cite[Theorem 5.1]{torrascasas2024} to the pair $Z\subseteq Z$ and $X \subseteq Z$.
As proposed in~\cite{representative}, $\varepsilon$-representative subsets of a given dataset can be used to train neural networks guaranteeing an accuracy similar to (depending on $\varepsilon$) that of the original dataset.

\subsection{Hausdorff distance bounds from matching diagrams}

In this subsection, we obtain a couple of bounds to the Hausdorff distance $d^{\scst Z}_{\scst H}(X, Z)$ between a subset $X$ and the metric space $Z$. 
The first, which is the weaker, is solely based on $\cMXZ$ and the infinity line of $\Df$ with the Hausdorff distance.
The second uses the computed information from $\TMT(Z)$ to deduce a sharper bound obtained with the computation of $\cMXZ$ described in this work. 
The results provided in this subsection are proven in~\ref{sec:proof-sharper-bound}.

\begin{proposition}\label{prop:bound}
    $\eta_f \leq d^{\scst Z}_{\scst H}(X, Z) \leq \sum_{b>0} \mf((\infty,b))
    b$.
\end{proposition}

Next, we define a quantity that can be obtained using $\TMT(Z)$. 
Let $c_f$ denote the maximum number of points from $Z\setminus X$ contained in a single connected component from $\VR_{\eta_f}(Z)$.
In addition, assuming $\eta_f\neq 0$, we denote 
\[\mbox{$
b_f = \sup\,\big\{ b>0 \mbox{ such that } \sum_{r\geq b} \mf((\infty,r))
\geq c_f\big\}\;$
and also}\] 
\[\mbox{$r_f = \sum_{b\geq b_f} \mf((\infty,b)) - c_f$,}\]
noticing that $b_f>0$ and $0 \leq r_f < 
\mf((\infty,b_f))$.

\begin{proposition}\label{prop:sharper-bound}
    $d^{\scst Z}_{\scst H}(X, Z) \leq 
    \Big(\sum_{b \geq b_f}\mf((\infty,b)) b\Big) - r_fb_f
    $.
\end{proposition}

Intuitively, the upper bound from Proposition~\ref{prop:sharper-bound} is computed by adding the $c_f$ largest quantities along the right infinity line from $\Df$.


\section{Methodology and applications }\label{applications}

To introduce the general ML framework used in the following experiment subsections, we first provide a rough explanation of the needed ML models and training to state the methodology used. 
As a deep explanation of ML techniques is beyond the scope of this paper, we 
refer to~\cite{Goodfellow-et-al-2016} for further details.

\subsection{ML preliminaries}

In our case, we will focus on supervised classification based on a dataset $Z$ given by points in $\mathbb{R}^n$, where $n$ is the number of features, labeled into different classes. Then, from a subset $X$ of the dataset, which is called the training set, we expect to generalize and be able to classify the remaining part of the dataset, called the test set. There are many ML models, but only multilayer perceptrons (MLPs) are used here. 
Specifically, each layer of a MLP consists of several interconnected neurons, and each neuron in a layer performs a computation that consists of a linear combination of parameters called weights and an activation function. 
We denote a given MLP by a product of factors where the number of factors is the number of layers and each factor is the number of neurons of the corresponding layer.
The parameters of the MLP are then trained following a training algorithm. In this paper, Adam \cite{KingBa15} was used as a training algorithm and sparse categorical cross entropy as a loss function. Finally, the trained model is evaluated on the test set. 
Given the qualitative nature of our proposed technique, we used a confusion matrix (with rows representing predicted classes and columns representing expected classes) to evaluate the model's performance and provide insight into the obtained results.

\subsection{Experiments}

In this section, we apply the concepts introduced in this paper to two different datasets: the housing dataset~\cite{housing}
and the dry beans dataset~\cite{dry-beans}.
The housing dataset is a standard dataset in machine learning and statistics, which contains various attributes related to housing prices and features. 
In contrast, the dry beans dataset encompasses morphological characteristics of different dry bean varieties, serving as an example from agricultural research. Both of them are tabular datasets, and their features are real-valued.

\begin{figure}[ht!]
    \centering
    \includegraphics[width=0.24\textwidth]{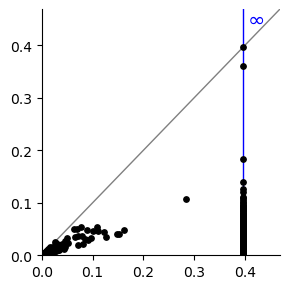}
    \qquad
    \includegraphics[width=0.24\textwidth]{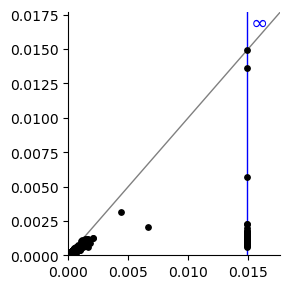}
    \caption{On the left (resp. on the right): 
    Depiction of the matching diagram $D({\cal H})$ (resp. $D({\cal B})$)
    associated to the 
    housing dataset $Z^{\cal H}$ 
    (reps. $Z^{\cal B}$) and a random subset $X^{\scst {\cal H}}$ (resp. $X^{\scst {\cal B}}$).
    The axes are scaled differently for each dataset.
    } \label{fig:housing_matching}
\end{figure}

Regarding the  housing dataset $Z^{\scst \cal H}$ and a random subset of it, $X^{\scst \cal H}$, let $D({\cal H})=D(X^{\scst \cal H},Z^{\scst \cal H})$. 
We  can see in Figure~\ref{fig:housing_matching} on the left, that most of the points of 
$D({\cal H})_{\scst O}$ are concentrated at $(0,0)$  (meaning that matched intervals are small) and that most of the points of $D({\cal H})_{\infty}$ at the point $(\infty,0)$ (meaning that unmatched intervals are small). 
Nevertheless, one matched interval difference reaches almost $0.3$, and one unmatched interval has a length between $0.3$ and $0.4$.
We can then expect some distortion of the original shape of the dataset with respect to the shape of the subset.
Regarding the  dry beans dataset $Z^{\cal B}$ and a random subset of it, $X^{\scst \cal B}$, in Figure~\ref{fig:housing_matching} on the right, $D({\cal H})$
shows that, in general, intervals of $\B(\PH_0(X^{\scst {\cal B}}))$ and $\B(\PH_0(Z^{\scst {\cal B}}))$ are well matched. 

Next, for both datasets, our approach involved training an MLP using a subset of the dataset and evaluating its performance on the remaining set as a test set to assess the generalization ability of the MLP. 
Then, for each class, we do a qualitative study based on the differences between the matched intervals. 
This gives us an idea of how well each of the classes is represented by the subset.
These experiments aim to show the robustness of the training process with respect to the training dataset, testing and highlighting the potential limitations of model performance evaluation when considering the dataset as a whole versus a subset.


\subsection{Housing dataset}

\begin{table}[ht]
\centering
\begin{adjustbox}{max width=0.4\textwidth}
\begin{tabular}{|c|c|c|c|c|}
\hline
Class & 0 & 1 & 2 & Total \\ 
\hline
0 & $\begin{array}{c}
       6994  \\
    {35.61\%}
  \end{array}$\cellcolor{blue!40} & $\begin{array}{c}
       1788  \\
    \textcolor{red}{9.10\%}
  \end{array}$\cellcolor{blue!20} & $\begin{array}{c}
       114  \\
    \textcolor{red}{0.58\%}
  \end{array}$\cellcolor{blue!10} & $\begin{array}{c}
       8896  \\
    \textcolor{green}{78.62\%} \\
    \textcolor{red}{21.38\%}
  \end{array}$\cellcolor{blue!60} \\ 
\hline
1 & $\begin{array}{c}
       2065  \\
    \textcolor{red}{10.51\%}
  \end{array}$\cellcolor{blue!20} & $\begin{array}{c}
       5798  \\
    {29.52\%}
  \end{array}$\cellcolor{blue!40} & $\begin{array}{c}
       1469  \\
    \textcolor{red}{7.48\%}
  \end{array}$\cellcolor{blue!30} & $\begin{array}{c}
       9332  \\
    \textcolor{green}{62.13\%} \\
    \textcolor{red}{37.87\%}
  \end{array}$\cellcolor{blue!50} \\ 
\hline
2 & $\begin{array}{c}
       14  \\
    \textcolor{red}{0.07\%}
  \end{array}$\cellcolor{blue!10} & $\begin{array}{c}
       237  \\
    \textcolor{red}{1.21\%}
  \end{array}$\cellcolor{blue!15} & $\begin{array}{c}
       1161  \\
    {5.91\%}
  \end{array}$\cellcolor{blue!25} & $\begin{array}{c}
       1412  \\
    \textcolor{green}{82.22\%} \\
    \textcolor{red}{17.78\%}
  \end{array}$\cellcolor{blue!45} \\ 
\hline
Total & $\begin{array}{c}
       9073  \\
    \textcolor{green}{77.09\%} \\
    \textcolor{red}{22.91\%}
  \end{array}$\cellcolor{blue!60} & $\begin{array}{c}
       7823  \\
    \textcolor{green}{74.11\%} \\
    \textcolor{red}{25.89\%}
  \end{array}$\cellcolor{blue!50} & $\begin{array}{c}
       2744  \\
    \textcolor{green}{42.31\%} \\
    \textcolor{red}{57.69\%}
  \end{array}$\cellcolor{blue!40} & $\begin{array}{c}
       19640  \\
    \textcolor{green}{71.04\%} \\
    \textcolor{red}{28.96\%}
  \end{array}$\cellcolor{yellow} \\ 
\hline
\end{tabular}
\end{adjustbox}
\caption{Housing dataset.
Confusion matrix of the MLP trained on a random subset $X^{\scst {\cal H}}\subset Z^{\scst \cal H}$ of size $1000$ evaluated on the test set.     }
\label{tab:housing_matching}
\end{table}

\begin{figure}[ht!]
    \centering
    \subfloat[Class 0]{\includegraphics[width=0.24\textwidth]{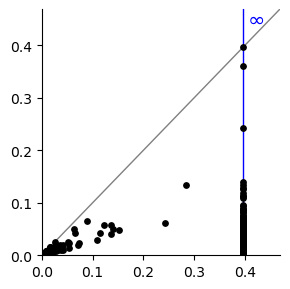}}
    \subfloat[Class 1]{\includegraphics[width=0.24\textwidth]{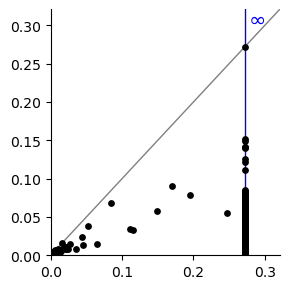}}    
    \subfloat[Class 2]{\includegraphics[width=0.24\textwidth]{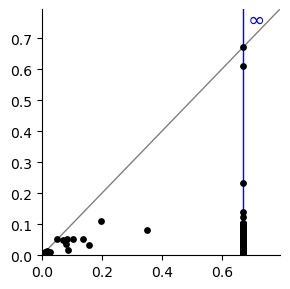}}
    \caption{Housing dataset. 
    Representation of the matching diagram 
    $D({\cal H}(i))$ for Class $i$ for $i\in [\![3]\!]$.
   The axes are scaled differently for each class.} 
    \label{fig:housing_matching_classes}
\end{figure}

The housing dataset  is a regression problem composed of 20640 real-valued 8-dimensional samples. We transformed it into a classification problem by discretizing the target values into a dataset $Z^{\cal H}$ of 3 classes and features normalized to $[0,1]$. 
We have chosen a random subset $X^{\scst {\cal H}}$ of size $1000$, and a $256\times 128 \times 64$ 
ReLu MLP with a final output layer with Softmax activation function, that was trained for $1000$ epochs using Adam training algorithm with a learning rate $0.001$. 
We repeated the training $10$ times, and the highest accuracy values we reached were $0.71$ on the test set and $0.72$ on the training set. For a more thorough study of the accuracy by classes for one of the iterations, see Table~\ref{tab:housing_matching}. 
The confusion matrix shows that Class 2 has less accuracy than the other classes. 

In Figure~\ref{fig:housing_matching_classes}, we provide the matching diagram for each class.
We denote by $D({\cal H}(i))$ the matching diagram corresponding to the restriction to the points of $X^{\scst {\cal H}}$ and $Z^{\scst \cal H}$ of Class $i$. 
As we can see, for class $2$, the points of
$D({\cal H}(2))_{\scst O}$
tend to disperse to the right and far from the diagonal axis, which means bigger differences between the matched intervals;
and one point of $D({\cal H}(2))_{\infty}$
(vertical line on the right) appear far from the $x$-axis.
In particular, a point appears at $y$-value $0.6$, which is much higher than any other point in any other class. 
Nevertheless, for class $0$, the points of
$D({\cal H}(0))_{\scst O}$ are more concentrated at $(0,0)$. 
This agrees with what we obtained in the confusion matrix of the trained MLP, represented in Figure~\ref{tab:housing_matching}.


\subsection{Dry beans dataset}

\begin{table}[ht!]
    \centering
    \resizebox{0.8\textwidth}{!}{
\begin{tabular}{|c|c|c|c|c|c|c|c|c|}
\hline
Class  & 0 & 1 & 2 & 3 & 4 & 5 & 6 & Total \\ 
\hline
0 & $\begin{array}{c}
       252 \\
       \textcolor{red}{2.37\%}
    \end{array}$\cellcolor{blue!20} &  & $\begin{array}{c}
       224 \\
       \textcolor{red}{2.11\%}
    \end{array}$\cellcolor{blue!20} &  & $\begin{array}{c}
       77 \\
       \textcolor{red}{0.73\%}
    \end{array}$\cellcolor{blue!10} &  &  & $\begin{array}{c}
       553 \\
       \textcolor{Green}{45.57\%} \\
       \textcolor{red}{54.43\%}
    \end{array}$\cellcolor{blue!35} \\ 
\hline
1 &  & $\begin{array}{c}
       405 \\
       \textcolor{red}{3.82\%}
    \end{array}$\cellcolor{blue!30} &  &  &  &  &  & $\begin{array}{c}
       405 \\
       \textcolor{Green}{100.00\%} \\
       \textcolor{red}{0.00\%}
    \end{array}$\cellcolor{blue!30} \\ 
\hline
2 & $\begin{array}{c}
       431 \\
       \textcolor{red}{4.06\%}
    \end{array}$\cellcolor{blue!30} & $\begin{array}{c}
       2 \\
       \textcolor{red}{0.02\%}
    \end{array}$\cellcolor{blue!5} & $\begin{array}{c}
       930 \\
       \textcolor{red}{8.76\%}
    \end{array}$\cellcolor{blue!40} &  & $\begin{array}{c}
       16 \\
       \textcolor{red}{0.15\%}
    \end{array}$\cellcolor{blue!5} &  &  & $\begin{array}{c}
       1379 \\
       \textcolor{Green}{67.44\%} \\
       \textcolor{red}{32.56\%}
    \end{array}$\cellcolor{blue!40} \\ 
\hline
3 &  &  &  & $\begin{array}{c}
       2619 \\
       \textcolor{red}{24.68\%}
    \end{array}$\cellcolor{blue!60} & $\begin{array}{c}
       78 \\
       \textcolor{red}{0.74\%}
    \end{array}$\cellcolor{blue!10} & $\begin{array}{c}
       934 \\
       \textcolor{red}{8.80\%}
    \end{array}$\cellcolor{blue!40} & $\begin{array}{c}
       328 \\
       \textcolor{red}{3.09\%}
    \end{array}$\cellcolor{blue!20} & $\begin{array}{c}
       3959 \\
       \textcolor{Green}{66.15\%} \\
       \textcolor{red}{33.85\%}
    \end{array}$\cellcolor{blue!60} \\ 
\hline
4 & $\begin{array}{c}
       289 \\
       \textcolor{red}{2.72\%}
    \end{array}$\cellcolor{blue!20} &  & $\begin{array}{c}
       107 \\
       \textcolor{red}{1.01\%}
    \end{array}$\cellcolor{blue!10} &  & $\begin{array}{c}
       719 \\
       \textcolor{red}{6.78\%}
    \end{array}$\cellcolor{blue!30} & $\begin{array}{c}
       22 \\
       \textcolor{red}{0.21\%}
    \end{array}$\cellcolor{blue!5} & $\begin{array}{c}
       98 \\
       \textcolor{red}{0.92\%}
    \end{array}$\cellcolor{blue!10} & $\begin{array}{c}
       1235 \\
       \textcolor{Green}{58.22\%} \\
       \textcolor{red}{41.78\%}
    \end{array}$\cellcolor{blue!35} \\ 
\hline
5 & $\begin{array}{c}
       3 \\
       \textcolor{red}{0.03\%}
    \end{array}$\cellcolor{blue!5} &  & & $\begin{array}{c}
       115 \\
       \textcolor{red}{1.08\%}
    \end{array}$\cellcolor{blue!10} & $\begin{array}{c}
       24 \\
       \textcolor{red}{0.23\%}
    \end{array}$\cellcolor{blue!5} & $\begin{array}{c}
       211 \\
       \textcolor{red}{1.99\%}
    \end{array}$\cellcolor{blue!20} & $\begin{array}{c}
       227 \\
       \textcolor{red}{2.14\%}
    \end{array}$\cellcolor{blue!20} &   $\begin{array}{c}
       580 \\
       \textcolor{Green}{36.38\%} \\
       \textcolor{red}{63.62\%}
    \end{array}$\cellcolor{blue!30} \\ 
\hline
6 & $\begin{array}{c}
       53 \\
       \textcolor{red}{0.50\%}
    \end{array}$\cellcolor{blue!10} &  & $\begin{array}{c}
       4 \\
       \textcolor{red}{0.04\%}
    \end{array}$\cellcolor{blue!5} & $\begin{array}{c}
       17 \\
       \textcolor{red}{0.16\%}
    \end{array}$\cellcolor{blue!5} &  $\begin{array}{c}
       604 \\
       \textcolor{red}{5.69\%}
    \end{array}$\cellcolor{blue!30} & $\begin{array}{c}
       402 \\
       \textcolor{red}{3.79\%}
    \end{array}$\cellcolor{blue!20} & $\begin{array}{c}
       1420 \\
       \textcolor{red}{13.38\%}
    \end{array}$\cellcolor{blue!35} & $\begin{array}{c}
       2500 \\
       \textcolor{Green}{56.80\%} \\
       \textcolor{red}{43.20\%}
    \end{array}$\cellcolor{blue!40}\\ 
\hline
Total & $\begin{array}{c}
       1028 \\
       \textcolor{Green}{24.51\%} \\
       \textcolor{red}{75.49\%}
    \end{array}$\cellcolor{blue!20} & $\begin{array}{c}
       407 \\
       \textcolor{Green}{99.51\%} \\
       \textcolor{red}{0.49\%}
    \end{array}$\cellcolor{blue!20} & $\begin{array}{c}
       1265 \\
       \textcolor{Green}{73.52\%} \\
       \textcolor{red}{26.48\%}
    \end{array}$\cellcolor{blue!30} & $\begin{array}{c}
       2751 \\
       \textcolor{Green}{95.20\%} \\
       \textcolor{red}{4.80\%}
    \end{array}$\cellcolor{blue!60} & $\begin{array}{c}
       1518 \\
       \textcolor{Green}{47.36\%} \\
       \textcolor{red}{52.64\%}
    \end{array}$\cellcolor{blue!20} & $\begin{array}{c}
       1569 \\
       \textcolor{Green}{13.45\%} \\
       \textcolor{red}{86.55\%}
    \end{array}$\cellcolor{blue!40} & $\begin{array}{c}
       2073 \\
       \textcolor{Green}{68.50\%} \\
       \textcolor{red}{31.50\%}
    \end{array}$\cellcolor{blue!20} & $\begin{array}{c}
       10611 \\
       \textcolor{Green}{61.78\%} \\
       \textcolor{red}{38.22\%}
    \end{array}$\cellcolor{yellow} \\ 
\hline
\end{tabular}
}  
    \caption{Dry beans dataset. Confusion matrix of the MLP trained on a random subset $X^{\scst {\cal B}}\subset {\cal B}$ of size $3000$ evaluated on the test set. }
    \label{tab:beans_conf_matrix}
\end{table}

\begin{figure}[ht!]
    \centering
\subfloat[Class 0]{\includegraphics[width=0.24\textwidth]{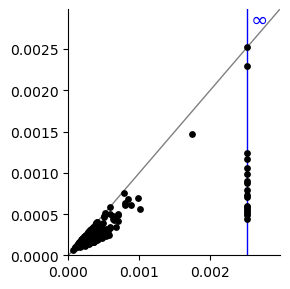}}
    \subfloat[Class 1]{\includegraphics[width=0.24\textwidth]{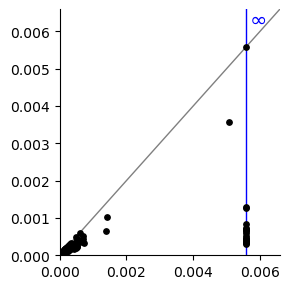}}
    \subfloat[Class 2]{\includegraphics[width=0.24\textwidth]{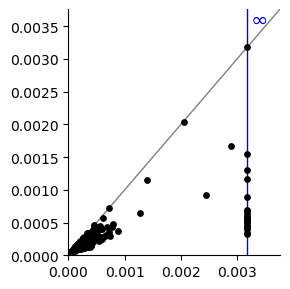}}
    \subfloat[Class 3]{\includegraphics[width=0.24\textwidth]{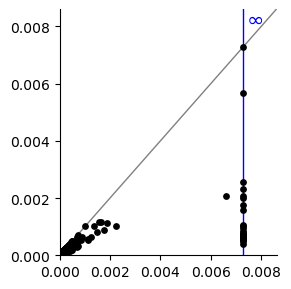}}\\
    \subfloat[Class 4]{\includegraphics[width=0.24\textwidth]{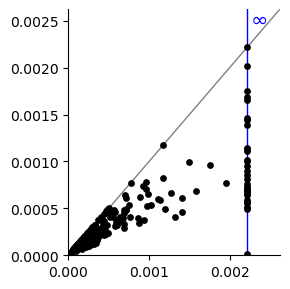}}
    \subfloat[Class 5]{\includegraphics[width=0.24\textwidth]{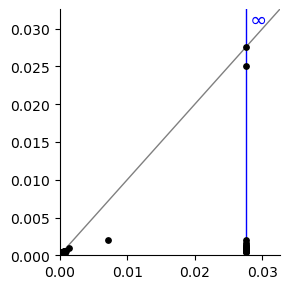}}
\subfloat[Class 6]{\includegraphics[width=0.24\textwidth]{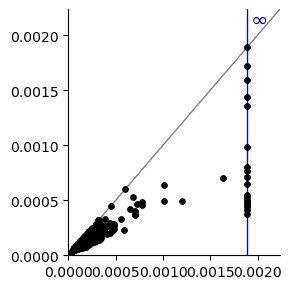}}
    \caption{Dry beans dataset. 
    Depiction of the matching diagram $D({\cal B}(i))$ for $i\in [\![7]\!]$
    corresponding to the restriction to the points of $X^{\scst {\cal B}}$
and $Z^{\scst \cal B}$ belonging to Class $i$.  
   }
    \label{fig:beans_matching_class}
\end{figure}

\begin{figure}[ht!]
    \centering
\begin{minipage}{.30\textwidth}
\begin{subfigure}{\textwidth}
\includegraphics[width=\textwidth]{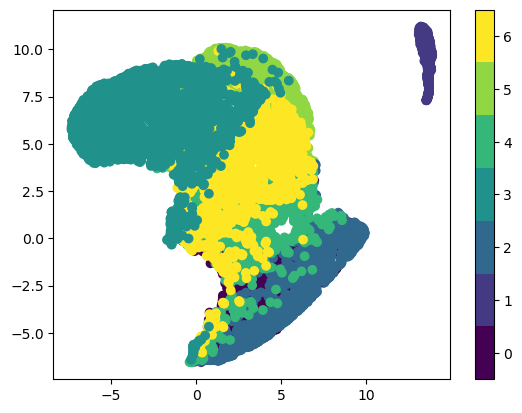} 
\end{subfigure}
\end{minipage}
   \begin{minipage}{.69\textwidth}
\begin{subfigure}{\textwidth}
\includegraphics[width=\textwidth]{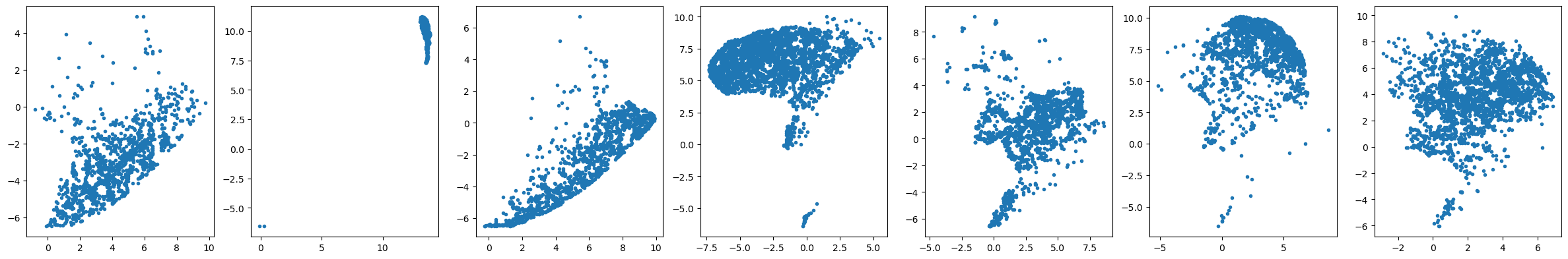} 
\end{subfigure}
\begin{subfigure}{\textwidth}
\includegraphics[width=\textwidth]{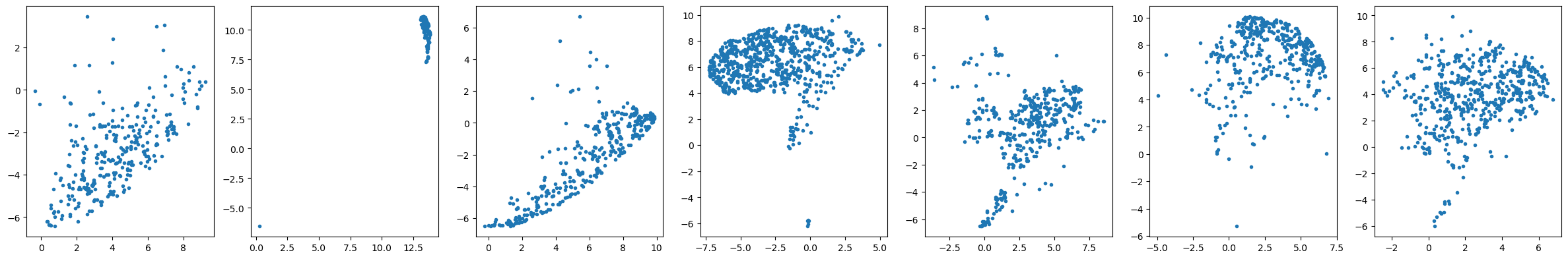}
\end{subfigure}
\end{minipage}
\hfill

    \caption{On the left, we show the UMAP embedding of the dataset colored by classes. On the right and from left to right, we show the embeddings for the different classes (top row) together with their subset (bottom row). 
    } 
    \label{fig:umap_embedding}
\end{figure}

The dry beans dataset $Z^{\scst \cal B}$ is composed of 13611 real-valued 16-dim. samples obtained from measures taken from 
images of $7$ different types of grains. 
In this experiment, we have  
chosen a random subset $X^{\scst {\cal B}}$ of size $3000$, leaving the remaining set as a test set, and trained a 
$1024\times 256\times 128 \times 64$ 
ReLu MLP with a final output layer with Softmax activation function.
The MLP was trained for $1000$ epochs using Adam training algorithm with learning rate $0.001$. We repeated the training $10$ times, and the highest accuracy value we reached was 0.63 on the test and training set. Table~\ref{tab:beans_conf_matrix} shows the confusion matrix on the test set for one of the iterations. 
There, we can see the accuracy for each of the classes. In this case, we can find that the class with the lowest performance is Class 5, and the one with the highest performance is Class 1. If we check the proportion of each class on the dataset, we appreciate that Class 1 (the class with the best performance) is the less represented, which seems contradictory. 

Let $D({\cal B})=D(X^{\scst \cal B},Z^{\scst \cal B})$.
In Figure~\ref{fig:beans_matching_class}, we can see that all the classes display similar matching diagrams and that points of $D({\cal B}(i))_{\scst O}$
for $i\in[\![6]\!]$ concentrate at $(0,0)$. 
Let us remark that the order of magnitude for the $x$-axis is $10^{-3}$ or lower for all classes except for Class 5. Also, we can see higher values for Class 5 on the $y$-axis and a point in $D({\cal B}(5))_{\infty}$
representing an unmatched interval of length close to $0.025$.
On the contrary, we can see that $D({\cal B}(1))_{\infty}$
presents all points close to $(\infty,0)$.
This could explain the poor (resp. good) performance of the trained MLP for Class 5 (resp. Class 1).

To check the intuition provided by 
$D({\cal B})$,
we have computed a 2-dimensional embedding of the dataset 
$Z^{\scst \cal B}$ (see Figure~\ref{fig:umap_embedding}) using the UMAP algorithm \cite{mcinnes2020umapuniformmanifoldapproximation}. 
The embedding shows how Class 1 is isolated by the embedding, making it easier to classify by the model. 


\section{A matrix computation procedure for $\cMXZ$}\label{sec:algorithm}

In this section, we adapt the matrix procedure explained in~\cite{matchings} to calculate the induced block function $\cMXZ:S^{\scst V}\times S^{\scst U}\rightarrow \Z_{\geq 0}$ previously defined. 


\subsection{Computational procedure overview}\label{subsec:overview}

Below, we provide an overview of the procedure to compute $\cMXZ$.
As in Subsection~\ref{subsec:interpretation-Df}, we index the points as $X=\{ z_0, \ldots , z_{\ell-1}\}$ and also $Z\setminus X = \{ z_{\ell}, \ldots , z_{n-1}\}$.
We compute the persistence modules   $V=\PH_0(X)$ and $U=\PH_0(Z)$, and the barcodes $\B(V)$ and $\B(U)$, using  minimum spanning trees of $\VR(X)$ and $\VR(Z)$
and 
the 
triplet merge trees obtained from them, $\TMT(X)$ and $\TMT(Z)$. 
In particular, we have that 
$V_0=\big\langle [z_0],\dots,[z_{\ell-1}]\big\rangle$ and 
$U_0=\big\langle [z_0],\dots,[z_{n-1}] \big\rangle$.

Next, we consider an order in $\TMT(X)$ where 
$(z_j, b_j, z_i)<(z_k, b_k, z_r)$ if and only if either $b_j < b_k$ or ($b_j=b_k$ and $z_j < z_k$) and fix the same order in $\TMT(Z)$.
Using Equation~(\ref{eq:decomposition-U-mUb}) and Remark~\ref{rem:iso-TMT-RepU},
we obtain a persistence isomorphism
\[\mbox{$\beta:
\left(\bigoplus_{
(z_j, b, z_i) \in \TMT(Z)}
\kk_{b} \right) \oplus \kk_\infty
\rightarrow U\,,$}
\]
defined for each component 
$(z_j, b, z_i)\in \TMT(Z)$,
by the morphism $\kk_{b} \rightarrow U$ that sends $1 \in 
(\kk_{b})_0$ to $[z_j]+[z_i] \in U_0$ 
together with $\kk_\infty \rightarrow U$ sending $1 \in (\kk_{\infty})_0$ to $[z_0]$.
Analogously, we consider the 
persistence isomorphism 
\[\mbox{$
\alpha: \left(\bigoplus_{(z_k, a, z_r)\in \TMT(X)} \kk_a\right) \oplus \kk_\infty
\rightarrow V\,.$}
\]
Thus, we can consider the composition 
\[\mbox{$
\beta^{-1}  f \alpha: 
\left( 
\bigoplus_{(z_k, a, z_r)\in \TMT(X)} \kk_a
\rightarrow 
\bigoplus_{(z_j, b, z_i)\in \TMT(Z)} \kk_b
\right) \oplus \Big( \kk_\infty\rightarrow \kk_\infty\Big)
\,.$}
\] 

We pay attention to the matrix $F$ associated with the linear map
$(\beta^{-1}  f \alpha)_0$ ignoring the component $(\kk_\infty)_0\rightarrow (\kk_\infty)_0$.
Since, by Remark~\ref{rem:iso-TMT-RepU}, 
we have that $\TMT(X) \simeq \Rep\B(V)$ and $\TMT(Z) \simeq \Rep\B(U)$ then
$F$ is a $\#\Rep\B(V)\times \#\Rep\B(U)$ matrix with columns indexed by $\TMT(X) \simeq \Rep\B(V)$ and rows indexed by $\TMT(Z) \simeq \Rep\B(U)$ and with coefficients in $\Z_2$. 
We call $F$ the \emph{associated matrix} to $f$.

Finally, we compute $R$, the Gaussian reduction of $F$, obtained by left-to-right column additions.
This leads to a computation for $\cMXZ(a,b)$, for $a \geq 0$ and $b \geq 0$, as 
\[
\cMXZ(a,b) = 
\# \left\{
\begin{array}{l}
\mbox{ pivots from } R \mbox{ in columns associated to } 
([0,a),t) 
\\
\mbox{ and rows associated to }
([0,b),r) 
\\ 
\mbox{ for all }
t \in \llbracket \mV\!(a)\rrbracket
\mbox{ and }
r\in \llbracket \mU\!(b)\rrbracket
\end{array}
\right\}\,.
\]

Summing up, we compute $\cMXZ$ with the following operations:
\begin{enumerate}
    \item\label{step:min-span-tree} Compute minimum spanning trees $\MST(X)$ and $\MST(Z)$ of $\VR(X)$ and $\VR(Z)$
    respectively.
    \item\label{step:tmt} Compute 
    the triplet merge trees $\TMT(X)$ and $\TMT(Z)$ from such minimum spanning trees.
    \item\label{step:compute-F} Compute $F$ using $\TMT(X)$ and $\TMT(Z)$.
    \item\label{step:gauss} Perform a Gaussian-column reduction of $F$.
\end{enumerate}
Notice that Steps~\ref{step:min-span-tree} and~\ref{step:gauss} can be performed by standard methods that are widely well-known, so we do not discuss them.  
Thus, we differ a detailed explanation of Steps~\ref{step:tmt} and~\ref{step:compute-F} to the following subsection. 
Now, we consider an example.

\begin{example}\label{ex:matrix-computation-Mf}
In this example, we check that the block function from Example~\ref{ex:matching-example} for the pair $X\subseteq Z$ coincides with the one obtained by the abovementioned matrix procedure.
\\
First, suppose we have computed $F$ the associated matrix  to $f: \PH(X)\rightarrow \PH(Z)$
(the computation of $F$ is done in Example~\ref{ex:F})
where the entries are sorted 
from lower to higher endpoint,
that is, the rows correspond to the endpoints $0.05, 0.05, 1,03, 1.11$ and $1.19$ of the intervals of  $\B(U)$ and the columns correspond to the endpoints 
$1.19$ and $2.24$ of the intervals of  $\B(V)$.
Then, we compute the  
reduced matrix $R$ 
obtained using left-to-right column additions:
\[
F = 
\left[ \begin{array}{cc}
0 & 0 \\
0 & 0 \\
0 & 0 \\
1 & 0 \\
1 & 1
\end{array}
\right]
\longrightarrow
R = 
\left[ \begin{array}{cc}
0 & 0 \\
0 & 0 \\
0 & 0 \\
1 & 1 \\
1 & 0
\end{array}
\right]
\]
In particular, notice that the reduced matrix $R$ has pivots in the last two rows that correspond to pairing the two intervals of $S^{\scst V}$ with the two longest intervals of $S^{\scst U}$, obtaining $\cMXZ(1.19, 1.19)=\cMXZ(2.24,1.11)=1$.
\end{example}


\subsection{Associated matrix $F$ computation}\label{subsec:matrix-F-computation}

In this subsection, we introduce a procedure to obtain the matrix $F$ associated to $(\beta^{-1}  f \alpha)_0$, 
ignoring the term $(\kk_\infty)_0\rightarrow (\kk_\infty)_0$.

As outlined in Subsection~\ref{subsec:overview}, to compute $F$, first we need to execute Step~\ref{step:min-span-tree} and compute the minimum spanning trees $\MST(X)$ and $\MST(Z)$ of $\VR(X)$ and $\VR(Z)$, which we assume is already done.

Let us now focus on Step~\ref{step:tmt} that computes triplet merge trees $\TMT(X)$ and $\TMT(Z)$.
First of all, notice that, if Step~\ref{step:min-span-tree} is performed using Kruskal's method~\cite{IntroAlgorithms2022}, one could inspect the union-find data structures to directly produce $\TMT(X)$ and $\TMT(Z)$. 
Unfortunately, usual implementations of minimum spanning trees do not provide access to such union-find data structures.
It is worth pointing out that Steps~\ref{step:min-span-tree} and~\ref{step:tmt} could also be directly obtained using~\cite{triplets}. 
However,  for the sake of flexibility,
we believe it is beneficial to compute the triplet merge trees from any previous computation of minimum spanning trees independently of the method used.
Next, we describe how to execute Step~\ref{step:tmt} in a lightweight and efficient manner, which is (yet) another adaptation of Kruskal's algorithm. 
This step allows implementing 
computations of triplet merge trees in a short Python script with few requirements.

We start by considering $\TMT(Z)$ to be an empty set.
Also, we consider a vector $C$ of length $n$.
The $i$-th coordinate of $C$ is denoted as $C[i]$.
Initially, $C=(0,1,\ldots, n-1)$.
Next, we iterate over increasing values $b>0$.
For each $b$, we set $E_b$ as the set of edges from $\MST(Z)$ of length $b$.
Then, while $E_b\neq \emptyset$, we do the following operations:

\begin{enumerate}[label=(\roman*)]
    \item \label{step:min-edge}
    take an edge $[z_i,z_j] \in E_b$ such that $\min\{C[i],C[j]\}$
    is as small as possible; if there is more than one such edge, pick any;
    \item \label{step:joinC} set $m\gets \min\big\{C[i],C[j]\big\}$ and 
    $M\gets \max\big\{C[i],C[j]\big\}$ (here notice that $m \neq M$ since otherwise $\MST(Z)$ would not be a tree);
    \item \label{step:no1} add the triplet $(z_M, b, z_m)$ to $\TMT(Z)$;
    \item \label{step:change-C} range over $k\in[\![n]\!]$ and if $C[k]=M$ then set $C[k]\gets m$; 
    \item \label{step:no2} remove $[z_i,z_j]$ from $E_b$.
\end{enumerate}
Then, we continue iterating over increasing values $b>0$.
Notice that, at the end of the iteration over a fixed $b>0$, for each point, $z_i \in Z$, the entry $C[i]$ is equal to the minimum index in $A_b(z_i)$; 
were $A_b(z_i)$ is the component in $\VR_b(Z)$ that contains $z_i$.
We  obtain $\TMT(X)$ by performing the same procedure on 
$\MST(X)$. 
The complexity of this step is 
${\cal O}\big (n\log n\big)$ since it is an adaptation of Kruskal's algorithm and $\MST(Z)$ has $n-1$ edges being $n=\# Z$.

Here, we obtain a result that justifies that our procedure leads to $\TMT(Z)$. 
In addition, we need this result to prove
Proposition~\ref{prop:sharper-bound}
in~\ref{sec:proof-sharper-bound}. 

\begin{proposition}\label{prop:edge-triplet-bijection}
    Denote by $\MST_b(Z)$ the subtree of $\MST(Z)$ containing the edges of length $\leq b$ and denote by $E$
    the set of edges of $\MST(Z)$.
    Given $(z_j, b, z_i)\in \TMT(Z)$, we have that $z_i$ and $z_j$ lie in the same component in $\MST_b(Z)$, which is represented by $z_i$. 
    Furthermore, there exists a bijection 
    $\alpha:\TMT(Z)\rightarrow E$ such that, for a given $(z_M,b,z_m) \in \TMT(Z)$, the edge $\alpha((z_M,b,z_m))$ has length $b$ and $z_m$ and $z_M$ are not path connected in 
    $\MST(Z)\setminus \{\alpha((z_M,b,z_m))\}$.  
\end{proposition}

\begin{proof}
    We consider the algorithm for computing $\TMT(Z)$ from $\MST(Z)$ outlined above and consider the iteration step at value $b>0$. 
    Recall that we iterate over edges $[z_i, z_j]$ 
    in $E_b$ adding triplets $(z_M, b, z_m)$ to $\TMT(Z)$ by Steps~\ref{step:joinC} and~\ref{step:no1}.
    Now, by Steps~\ref{step:min-edge} and~\ref{step:change-C}, the value $m$ 
    cannot decrease as we iterate over $E_b$. 
    This implies that  $z_M$ and $z_m$ lie in the same component in $\MST_b(Z)$ represented by $z_m$.
\\
    Next, we show the second claim. Consider again the iteration step at value $b$, where we range over edges 
    $[z_i, z_j]$ 
    in $E_b$ adding triplets $(z_M, b, z_m)$ to $\TMT(Z)$. 
    This implies that there exists a path $\gamma$ in $\MST_b(Z)$ joining $z_M$ and $z_m$ and having 
    $[z_i, z_j]$ 
    as an edge. 
    We set $\alpha((z_M, b, z_m))=[z_i, z_j]$.
    Now, 
    $z_i$ and $z_j$
    are not path connected in $\MST(Z)\setminus \{\alpha((z_M, b, z_m))\}$ since otherwise $\MST(Z)$ could not be a tree. 
        \qed
\end{proof}

Now we focus on Step~\ref{step:compute-F} that computes the matrix $F$ associated to 
$(\beta^{-1}  f \alpha)_0$.
 Notice that one could compute $F$ from Step~\ref{step:min-span-tree} by using the first differential matrices from $\MST(X)$ and $\MST(Z)$. However, the advantage of using triplet merge trees is that it leads to a faster Gaussian reduction because the matrices involved have fewer rows.

We first consider the canonical basis for $V_0=\Ho_0(\VR_0(X))$ where generators correspond to points from $X$. 
We denote by $A_0$ the matrix associated with $\alpha_0$ ignoring the summand $(\kk_\infty)_0$. 
That is, $A_0$ has rows indexed by 
$X$ and columns indexed by $\TMT(X)$, where the $(z_j,b,z_i)$-column has non-trivial entries in the positions $i$ and $j$.

Next, we consider the canonical basis for $U_0=\Ho_0(\VR_0(Z))$, where generators are given by points from $Z$.
Further, $U_0$ can be decomposed into a direct sum  $U_0 = U_0^{\scst X} \oplus U_0^{\scst Z\setminus X}$ where $U_0^{\scst X}$ is generated by the points from $X$ and $U_0^{\scst Z\setminus X}$ is generated by the points from $Z\setminus X$. 
Now, consider the following linear map,  
\[\mbox{$
\beta^{\scst X}_0: 
\left(\bigoplus_{(z_j,b,z_i)\in \TMT(X,Z)} 
(\kk_{b})_0 \right) \oplus (\kk_\infty)_0
\rightarrow 
U^{\scst X}_0$.}
\]
which is a restriction of $\beta_0$, 
and so it is injective. 
Further, $\beta^{\scst X}_0$ is a bijection since $\#\TMT(X,Z)=\# X - 1$.
Now, we denote by $B_0$ the associated matrix to $\beta^{\scst X}_0$ ignoring the $(\kk_\infty)_0$ component.
That is, $B_0$ has rows indexed by $X$ and columns indexed by $\TMT(X, Z)$, where the $(z_j, b, z_i)$-column has non-trivial entries in the positions $i$ and $j$.

Next, as $f_0$ is induced by the inclusion $X\subseteq Z$, we have that $f_0$ is equal to $f^{\scst X}_0:V_0 \rightarrow U_0^{\scst X}$ composed with the inclusion $\iota^{\scst X}_{\scst U}: U_0^{\scst X}\hookrightarrow U_0$. 
We denote by $F^{\scst X}$ the matrix associated to $(\beta^{\scst X}_0)^{-1}\circ f_0^{\scst X} \circ \alpha_0$, ignoring the terms $(\kk_\infty)_0$.
Adding zero rows corresponding to triples in $\TMT(Z)\setminus \TMT(X,Z)$ to $F^{\scst X}$, we obtain $F$. To see why, notice that there is a commutative diagram 
\[
\begin{tikzcd}[column sep=0.55cm
]
&
\left(\bigoplus_{\scst (z_j, b, z_i)\in \TMT(X)} (\kk_b)_0 \right) \oplus (\kk_\infty)_0 
\ar[ld, "\alpha_0", swap]
& 
\\
V_0 \ar[r, "f_0"] \ar[rd, "f^{\scst X}_0", swap
]
& 
U_0 \ar[r, "(\beta_0)^{-1}"]
& 
\left(\bigoplus_{\scst (z_j, b, z_i)\in \TMT(Z)} (\kk_b)_0\right) \oplus (\kk_\infty)_0
\\ 
& 
U^{\scst X}_0 \ar[r, "(\beta^{\scst X}_0)^{-1}"]
\ar[u, "\iota^X_U", swap]
& 
\left(\bigoplus_{\scst (z_j, b, z_i)\in \TMT(X,Z)} (\kk_b)_0\right) \oplus (\kk_\infty)_0
\ar[u, "\iota_{\TMT}^X", swap]
\end{tikzcd}
\]
where $\iota_{\scst \TMT}^{\scst X}$ is induced by the inclusion $\TMT(X,Z)\subseteq \TMT(Z)$.
In particular, we can rewrite the formula of $\cMXZ$ from Subsection~\ref{subsec:overview} as follows 
\begin{equation}
\label{def:Mf-triplets}
\cMXZ(a,b) = 
\# \left\{
\begin{array}{l}
\mbox{ pivots from } R \mbox{ in columns associated to } 
(z_i, a, z_j) 
\\
\mbox{ and rows associated to }
(z_s, b, z_t) 
\mbox{ for all }
\\ 
(z_i, a, z_j) \in \TMT(X)
\mbox{ and }
(z_s, b, z_t) \in \TMT(X,Z)
\end{array}
\right\}
\end{equation}

\begin{proposition}\label{prop:Mf-bijection-TMTX-TMTXZ}
    $\cMXZ$ induces a bijection $\TMT(X)\simeq \TMT(X,Z)$.
\end{proposition}

\begin{proof}
    We use Equation~(\ref{def:Mf-triplets}) and that $\# X -1= \# \TMT(X)=\#\TMT(X,Z)$. 
    Now, $F^{\scst X}$ is a square matrix 
    with columns indexed by $\TMT(X)$ and rows indexed by $\TMT(X,Z)$
    of size $\# X-1$ which is invertible, since it is the matrix associated to the composition of three isomorphisms $(\beta^{\scst X}_0)^{-1}\circ f_0^{\scst X} \circ \alpha_0$, concluding that  $\cM^{\scst Z}_{\scst X}$ induces a bijection.
    \qed
\end{proof}

Now, we explain how to compute $F^{\scst X}$.
We consider the block matrix below, which we reduce via left-to-right column additions:
\[
\left(\begin{array}{cc}
\id     & {\bf 0}         \\
B_0   &   A_0
\end{array}
\right)
\mapsto 
\left(\begin{array}{cc}
\id     & F^{\scst X}         \\
B_0   &   {\bf 0}
\end{array}
\right)
\]
where $\id$ is the identity matrix 
and ${\bf 0}$ is the zero matrix, both of dimension $(\ell-1)\times (\ell-1)$.
Notice that the pivots from $B_0$ are all unique, which allows a quick search for columns corresponding to pivots. 
The result obtained from the upper right corner of the reduced block matrix corresponds to the first $\ell-1$ rows, and the last $\ell-1$ columns equals $F^{\scst X}$. 
The worst-case complexity of this reduction is about $\cO(\ell^3)$, where $\ell=\# X$, since we have to reduce the matrix $A_0$ which 
has dimension  $\ell\times(\ell-1)$,
by adding the columns from $B_0$, which has already been reduced.
In practice,  the computation time is 
shorter than theoretically expected
since the matrices we deal with are sparse.

\begin{example}\label{ex:F}
    In this example, we explain how to compute the matrix $F$ that we used in Example~\ref{ex:matrix-computation-Mf}. 
    Recall the descriptions of $\TMT(X)$ and $\TMT(Z)$ 
    from Examples~\ref{ex:pershom-merge-trees} and~\ref{ex:barcode-PH0-Z}.
     Then, we might write $B_0$, $A_0$, and $F^{\scst X}$ as follows
    \[
    B_0 = \left[\begin{array}{ccccc}
        0 & 1 \\
        1 & 1 \\
        1 & 0 
    \end{array}\right]\,,\quad
    A_0 = \left[\begin{array}{cc}
        1 & 1 \\
        0 & 1 \\
        1 & 0 
    \end{array}\right]
    \mbox{ and }
    F^{\scst X}= 
    \left[\begin{array}{ccccc}
        1 & 0 \\
        1 & 1 
    \end{array}\right]
    \]
    where $B_0$ has columns indexed by 
   $ \TMT(X, Z)=\{(z_2, 1.11, z_1), (z_1, 1.19, z_0)\}$, which correspond to the last two triplets from $\TMT(Z)$, and rows indexed by $X=\{z_0,z_1,z_2\}$,   and $A_0$ has columns indexes by 
   $\TMT(X)=\{(z_2, 1.19,
   z_0)$, 
   $(z_1, 2.24, z_0)\}$ and rows indexes by $X$.
    Then, $F$ is obtained by adding $3=(\# Z)-(\# X)$ trivial rows on top of $F^{\scst X}$, as can be checked in Example~\ref{ex:matrix-computation-Mf}.
\end{example}

\begin{table}[ht!]
\centering
\begin{adjustbox}{max width=0.6\textwidth}
\begin{tabular}{ccccccc}
        \toprule
                & \textbf{Dimension}&  100 & 200 & 500 & 1000 \\
                        \midrule
        \textbf{Size} & \textbf{Proportion} & \textbf{ Time } & \textbf{ Time } & \textbf{ Time } & \textbf{ Time } 
        \\
        \midrule
        \multirow{4}{*}{1000}  & 0.1 & 0.2137 & 0.2251 & 0.2906 & 0.4167 \\
          & 0.2 & 0.2205 & 0.2337 & 0.3025 & 0.4266 \\
          & 0.5 & 0.2614 & 0.2770 & 0.3672 & 0.5181 \\
          & 0.8 & 0.3420 & 0.3647 & 0.4785 & 0.6791 \\\hline
        \multirow{4}{*}{5000}  & 0.1 & 6.8683 & 7.4391 & 9.3776 & 13.2139 \\
          & 0.2 & 7.1100 & 7.6029 & 9.6264 & 13.5654 \\
          & 0.5 & 8.3384 & 9.0599 & 11.3875 & 15.9382 \\
          & 0.8 & 11.0920 & 11.9642 & 15.1606 & 21.0903 \\\hline
        \multirow{4}{*}{10000} & 0.1 & 30.9312 & 33.8928 & 43.3496 & 59.5236 \\
         & 0.2 & 31.7921 & 34.5644 & 44.4252 & 61.1580 \\
         & 0.5 & 37.6867 & 41.3084 & 52.5283 & 72.4249 \\
         & 0.8 & 50.1831 & 54.3484 & 69.9975 & 96.1254 \\
        \bottomrule
    \end{tabular}
    \end{adjustbox}
\caption{Random datasets for different dimensions and sizes;  random subsets taking different proportions of the datasets. 
The values given are the execution time in seconds for the computation of $F$.
}
\label{tab:time_execution}
\end{table}

Now, we give an estimate of the complexity of our procedure. We have already discussed the complexity of Steps~\ref{step:tmt} and~\ref{step:compute-F}, so we first discuss the complexity of the remaining steps. For Step~\ref{step:min-span-tree}, since the number of edges in $\VR(Z)$ is about $n^2$, where $n=\# Z$, we can use the standard complexity of computing minimum spanning trees to obtain a complexity of $\cO(n^2\log n)$.
Now, for Step~\ref{step:gauss}, the worst case complexity for the Gaussian elimination of $F$ is $\cO(\ell^3)$, where $\ell=\# X$, since $F$
has  $n-\ell$ trivial rows.
Overall, the worst-time complexity of the four terms in the procedure is
$
\cO(n^2\log n + n\log n + \ell^3 + \ell^3)$ 
which simplifies to $\cO(n^2 \log n + \ell^3)$. 
In our experiments, the computation of Step~\ref{step:min-span-tree} takes the most time since usually, $\# Z$ is much larger than $\# X$. For example, in Table~\ref{tab:time_execution}, random sets of different sizes and dimensions were generated, and the time of execution of the proposed algorithm is shown.
Furthermore, these
times could be optimized in the future by considering the particular structures of the matrices we deal with.


\section{Conclusion and future works}
\label{sec:future}

In this paper, we have introduced a novel indicator of data quality. It is based on topological features
 obtained from 0-dimensional persistence modules and induced block functions. The goal of this tool is 
 to measure the quality of a training subset relative to the entire input dataset.
    The experimentation reveals that the topological data quality information is coherent with related tools to visualize/inspect datasets and demonstrates that it is a valuable measure, providing insights into why a chosen training subset might lead to poor performance.

Besides, we have demonstrated in~\cite{torrascasas2024} that 
the proposed definition of matching diagrams is stable.
We have not included here the proof of the stability property of $D(X,Z)$ for being rather technical and out of the scope of this paper.
Intuitively, the stability of $D(X,Z)$ implies that for a pair of samples $X \subseteq Z$, slight perturbations resulting in a new pair $X' \subseteq Z'$ will lead to approximately equivalent  matching diagrams  $D(X,Z)$ and $D(X',Z')$.

There are different research lines for future work. 
An interesting direction is to adapt our method to be robust to outliers and to optimize our code for large datasets. 
While our method is polynomial in theory, persistent homology is linear in practice, and we believe our method can achieve linear performance in practice too. 
Specifically, 
the current computational bottleneck of our implementation is in computing the minimum spanning tree in Step~\ref{step:min-span-tree}.
A future direction is to adapt the work from~\cite{triplets} to compute Steps~\ref{step:min-span-tree} and~\ref{step:tmt} directly, which is likely more efficient and parallelizable. 
Additionally, we have observed that the distance metric used to compute persistent homology should be related to the architecture chosen. 
For example, if we use Euclidean distance, it makes sense to use perceptrons.
If using the Structural Similarity Index Measure  or the Fréchet Inception Distance, a convolutional neural network should be used as the architecture. 
Moreover, in future work, we will also explore scenarios where $X$ is not a subset of $Z$ and topological features of dimensions higher than $0$ that can help in cases where 0-dimensional topological features are not decisive. 

\vspace{0.2cm}

\noindent{\bf Code Availability.}
The source code for the examples and experiments is available on the GitHub repository \cite{tdqual}.

\vspace{0.2cm}

\noindent{\bf Acknowledgments.}
Partially funded by the European Union under grant agreement no. 101070028-2
REXASI-PRO~\cite{rexasipro}, and by
MCIN/AEI and the NextGenerationEU/PRTR, under project TED2021-129438B-I00.


\bibliographystyle{plainnat}
\bibliography{biblio2}

\appendix 

\section{
Proofs}
\label{sec:proof-sharper-bound}

\subsection{Proof of Proposition~
\ref{prop:bound} }

Let us prove that $\mbox{$\eta_f \leq d^{\scst Z}_{\scst H}(X, Z) \leq \sum_{b>0} \mf((\infty,b))
    b$}.$
   First, observe that the upper bound is a consequence of Proposition~\ref{prop:sharper-bound}, where a sharper bound is given.
    Now, we prove the lower bound. For all $r < \eta_f$, by Proposition~\ref{prop:components-infty-line},
    there exists at least one component from $\VR_r(Z)$ with no samples from $X$. In particular, there exists at least one point $z \in Z\setminus X$ such that 
    $d^{\scst Z}(z, X)>r$. 
    Hence $\eta_f \leq 
    d^{\scst Z}(z, X)$ and the lower bound follows, since $d^{\scst Z}(z, X) \leq d^{\scst Z}_H(X,Z)$.

\subsection{Proof of Proposition \ref{prop:sharper-bound}}

We want to prove that $d^{\scst Z}_{\scst H}(X, Z) \leq 
    \Big(\sum_{b \geq b_f}\mf((\infty,b)) b\Big) - r_fb_f$
    .
    Suppose we have computed $\MST(Z)$ and $\TMT(Z)$ as directed in Subsection~\ref{subsec:overview}. 
     Let $\alpha\colon \TMT(Z) \to E$ be the bijection from Proposition~\ref{prop:edge-triplet-bijection}.
    Now, consider a point $z_i \in Z\setminus X$. 
    By definition, there exists a path (with no cycles) $\gamma$ in $\MST(Z)$ connecting $z_i$ with a point in $X$ and such that all edges have length $\leq \eta_f$. 
    Now, if it exists, we take any edge $e$ from $\gamma$ such that $e = \alpha((z_j, b, z_k))$ for $(z_j, b, z_k) \in \TMT(X,Z)$.
    If such $e$ does not exist, then fix $\gamma$. Otherwise, we modify $\gamma$ as follows. 
    By Proposition~\ref{prop:edge-triplet-bijection}, there exists a path $\tau$ in $\MST(Z)$ such that its edges have length $b\leq \eta_f$, it starts at $z_j$ and ends at $z_k$, and it goes along the edge $e$.  
    We combine $\gamma$ with $\tau$ to obtain a new path $\gamma'$ which joins $z_i$ with either $z_j$ or $z_k$ in such a way that $e$ is not in $\gamma'$. 
    Now, we look again for any edge (if it exists) $e'$ in $\MST(Z)$ such that $e' = \alpha((z_n, b, z_m))$ for $(z_n, b, z_m) \in \TMT(X,Z)$ and modify again $\gamma'$ if necessary. 
    Now, as we iterate, we cannot consider the same edge from $\alpha(\TMT(X,Z))$ twice since otherwise $\gamma'$ would have a cycle. 
    Since the number of triplets is finite, eventually, we must obtain a path $\gamma'$ joining $z_i$ to a point from $X$ and with no edges from $\alpha(\TMT(X,Z))$.
    Now, for each $z_i \in Z\setminus X$, we consider the resulting path $\gamma'$ joining $z_i$ to a point $z_j \in X$, and we assume that $\gamma'\cap X=\{z_j\}$, as otherwise $\gamma'$ can be shortened. 
      By the triangle inequality, we obtain 
$\dZ(z_i,z_j) \leq \sum_{(z_k, z_l) \in \gamma'} \dZ(z_k, z_l)$.
    Now, $\gamma'$ lies in a connected component from $\VR_{\eta_f}(Z)$, a path with maximum length $c_f$. Also, the lengths of the edges from $\gamma'$ injectively correspond to lengths of edges from $\alpha(\TMT(Z)\setminus \TMT(X,Z))$.

\end{document}